\documentclass[12pt,reqno]{amsart}

\usepackage[mathlines]{lineno}

\usepackage[english]{babel}
\usepackage{graphicx}
\usepackage{textcomp}
\usepackage{amssymb}
\usepackage{amsmath}
\usepackage{setspace}
\usepackage{geometry}
\usepackage{bm}
\usepackage{centernot}
\usepackage{tikz-cd}
\usepackage{hyperref}

\newtheorem{proposition}{Proposition}
\newtheorem{corollary}{Corollary}
\newtheorem{remark}{Remark}

\newtheorem{lemma}{Lemma}

\newcommand{\id}{{\mathbf 1}}
\newcommand{\supp}{\operatorname{supp}}

\usepackage{xcolor}

\begin{document}
\title{Mixed-Fourier-norm spaces and holomorphic functions}

%\author{Zhirayr Avetisyan \and Alexey Karapetyants}

\author{Zhirayr Avetisyan}
\address{Department of Mathematics: Analysis, Logic and Discrete Mathematics, Ghent University, 9000 Ghent, Belgium}
\email{jirayrag@gmail.com}

\author{Alexey Karapetyants}
\address{Institute of Mathematics, Mechanics and Computer Sciences $\&$ Regional Mathematical Center, Southern Federal
University, Rostov-on-Don, 344090, Russia}
\email{karapetyants@gmail.com}


\begin{abstract}
We describe a general framework of functional and Fourier analysis on domains with a free action of an Abelian Lie group $G$. Namely, on a domain of the form $G\times Y$ we introduce the appropriate spaces of distributions and measurable functions, establishing their most basic properties. Then we consider the half-Fourier transform $f(x,y)\mapsto\hat f(\xi,y)$ in the first variable, and discuss the behaviour of function spaces on $G\times Y$ and $\hat G\times Y$ under this transform. We introduce general mixed-Fourier-norm spaces on $G\times Y$, and the subspaces of holomorphic functions among them, and give an explicit descriptions of the Fourier images of these spaces.
\end{abstract}

\maketitle

%\noindent {\bf Keywords:}\,\, {good keywords}

%\noindent {\bf AMS MSC 2020:}\,\, 47B90 Operator theory and harmonic analysis ;

\section*{Introduction}

The spaces of functions defined in terms of conditions on Fourier images deserve special attention in analysis. This idea is not new and its implementation is dictated by considerations of the study of operators for which, in one sense or another, a realisation in terms of Fourier multipliers, or in terms close to this, is available.

Following this idea, in the series of works
\cite{KS-2018,KS-2017,KS-2016-1,KS-2016-2,KSmi-CVEE,KSmi-IzvVuz,
K-Grand-2021}, the spaces of holomorphic functions with conditions on the Fourier coefficients were investigated. These works are related to elliptic geometry in a disc and, accordingly, the Fourier coefficients are functions of a radial variable.

However, in the unit disc there are two more types of coordinate geometry - parabolic and hyperbolic, which are easier to visualise when carried over to the half-plane model, where they turn into Cartesian and polar coordinates on the half-plane, respectively. It would be natural to expect the development of a function space theory similar to the above works in these two cases.

The peculiarity of the elliptic geometry and Fourier coefficients in a disc is that in this case, holomorphic functions are uniquely determined by the Taylor series (Taylor coefficients), and this can be easily recalculated in terms of the Fourier coefficients. Such a criterion allows us to uniquely identify the subclass of holomorphic functions within the new spaces and constructively characterize these functions.

The analogous situation does not take place in the other two cases; theorems of the Paley-Wiener type can provide sufficient conditions, but do not provide a complete and constructive characterization. It is possible to operate with the Cauchy-Riemann conditions, which is actually done in this paper, but this approach requires a general view of the entire problem, which led to the studies that constitute the main result of this article.

Namely, a general approach to the definition and characterization of functions defined in terms of conditions on Fourier images is implemented, which covers all three geometries of a disc.

Let us take a moment to go back to earlier studies that also served as a motivation for this work.
In a series of papers \cite{Grudsky-Karapetyants-Vasilevski-JOpTh-2003, Grudsky-Karapetyants-Vasilevski-IEOT-2003, Vasilevski-Grudski-Karapetyants-SEMR-2006, Grudsky-Karapetyants-Vasilevski-BolSocMathMex-2004, Grudsky-Karapetyants-Vasilevski-JOpTh-2004}, a general theory of some classes of Toeplitz operators with special (generally unbounded) symbols in (weighted) Hilbert spaces of holomorphic functions (of the Bergman type) associated with three types of hyperbolic geometry in the unit disc was constructed. The developed methods and approaches work effectively in non-standard situations - in cases of non-standard growth of a function from the space under study, or ``bad'' behavior of the Toeplitz operator symbol (unboundedness, oscillation when approaching the boundary). These approaches are based on the study of the structural properties of the corresponding function spaces and on a special characterization of functions from these spaces. The authors used the transition to Fourier images for a constructive description of the subclass of holomorphic functions, as well as for characterizing the spectrum of the corresponding Toeplitz operator in each of the three cases. Let us repeat that this works effectively and exclusively in the case of a Hilbert space, when the transition to Fourier images is justified by the Parseval equation. This theory and some further aspects are best reflected in the book \cite{Vbook}.

In fact, this paper is in particular also an extension of the above mentioned scheme to the case where it is impossible to go directly to Fourier images, as was done in the case of a Hilbert space. Instead, we define our new spaces in Fourier images directly, which undoubtedly requires serious justification.

The structure of the paper is as follows. In Section 1 we lay out the geometrical setting: a product $G\times Y$, where $G$ is an Abelian Lie group and $Y$ a manifold (both assumed to be connected for the most part). The basis of the distribution theory on $G\times Y$ is the test function space
$$
\mathcal{D}(G\times Y)=\mathcal{S}(G)\,\hat\otimes\, C^\infty_c(Y),
$$
which consists of smooth functions that are of rapid decay (Schwartz) in $x\in G$ and of compact support in $y\in Y$. Thus, this space is not a standard one but a mixture; while its general properties are well-understood thanks to the abstract theory, here we introduce an explicit system of seminorms convenient for further considerations. The dual space $\mathcal{D}(G\times Y)'$ consists of distributions that are tempered in $x\in G$.

The central object of our study is the half-Fourier transform
$$
\mathfrak{F}:\mathcal{D}(G\times Y)\to\mathcal{D}(\hat G\times Y),
$$
or the Fourier transform in the variable $x$ only. Here $\hat G$ is the dual of $G$ - another Abelian Lie group, so that the usual Fourier transform $\mathfrak{F}$ maps $\mathcal{S}(G)$ isomorphically to $\mathcal{S}(\hat G)$. The half-Fourier transform $\mathfrak{F}:\mathcal{D}(G\times Y)\to\mathcal{D}(\hat G\times Y)$ is then the extension of the componentwise operator
$$
\mathfrak{F}\otimes\mathbf{1}:\mathcal{S}(G)\otimes C^\infty_c(Y)\to\mathcal{S}(\hat G)\otimes C^\infty_c(Y).
$$
The half-Fourier transform is further extended to
$$
\mathfrak{F}:\mathcal{D}(G\times Y)'\to\mathcal{D}(\hat G\times Y)'
$$
by duality, as usual. In the bulk of the paper, we will be concerned only with the half-Fourier transform on $G\times Y$, and will refer to it simply as the Fourier transform, without the risk of confusion.

In Section 2, we introduce the space $\mathcal{O}_M(G\times Y)$ of multipliers of $\mathcal{D}(G\times Y)$, i.e., the appropriate algebra of functions acting continuously on test functions by pointwise multiplication. Again, since this is a mixture of standard multiplier spaces, the derivations of basic properties are carried out in details. Here we also introduce the algebra $\mathrm{D}_G(G\times Y)$ of $G$-invariant differential operators, and the algebra $\mathrm{D}_M(G\times Y)$ of differential operators with coefficients growing moderately in the variable $x\in G$.

Section 3 is devoted to measurable functions on $G\times Y$ and Lebesgue spaces $L^p(G\times Y)$ with respect to a smooth $G$-invariant measure on $G\times Y$. The $L^2$ inner product $(\cdot,\cdot)_2$ gives rise to the operation of transpose on differential operators, $P(\partial)\mapsto P(\partial)^\top$, which is then used to define weak derivatives of distributions. It is important to note that, while the concrete identification of locally integrable functions with distributions, and the definitions of weak derivatives do depend on the choice of the measure, the results obtained in this paper using this identification and weak derivatives directly (Lemma \ref{CRHolLemma}, Proposition \ref{CRFourierProp}, Corollary \ref{HolDCorr} and Lemma \ref{hatu_0L1locLemma}) remain true regardless of the chosen measure, as long as the latter is given by an orientation (non-vanishing smooth top form).

Section 4 is where mixed-Fourier-norm spaces appear. For rather arbitrarily given Banach spaces of distributions $\mathcal{Y}(Y)$ on $Y$ and locally integrable functions $\Xi(\hat G)$ on $\hat G$, we first construct the space $\Xi(\hat G,\mathcal{Y}(Y))$ of distributions $\hat u(\xi,y)$ on $\hat G\times Y$, such that the map $\hat G\ni\xi\mapsto\hat u(\xi,\cdot)\in\mathcal{Y}(Y)$ is locally Bochner-integrable, and the map $\hat G\ni\xi\mapsto\|\hat u(\xi,\cdot)\|_{\mathcal{Y}(Y)}$ is in $\Xi(\hat G)$. We study the basic properties of this space under fairly mild assumptions. The corresponding mixed-Fourier-norm space $\mathcal{X}(G\times Y)$ is the Fourier-preimage of $\Xi(\hat G,\mathcal{Y}(Y))\cap\mathcal{D}(\hat G,\mathcal{Y}(Y))$, that is,
$$
\mathcal{X}(G\times Y)=\left\{u\in\mathcal{D}(G\times Y)'\;\vline\;\hat u=\mathfrak{F}u\in\Xi(\hat G,\mathcal{Y}(Y))\right\}.
$$
Note that in general this space is only a topological linear cone, but becomes a normed vector space under natural assumptions. The completeness of this space turns out to be subtle; here we will only show that it is closed in the topology of convergence in measure.

In Section 5 we introduce a complex structure on $G\times Y$ by assuming that $Y\subset\mathbb{R}^n$ is an open subset, $n=\dim G$, and viewing $G\times Y$ as $G+\imath Y$. This allows us to speak about holomorphic functions $u\in\mathrm{Hol}(G\times Y)$. The main result of this section is encoded in Proposition 4 and Corollary 4, where we show that a distribution $u\in\mathcal{D}(G\times Y)'$ is given by a holomorphic function $u\in\mathrm{Hol}(G\times Y)$ if and only if its Fourier transform factorizes as $\hat u(\xi,y)=e^{-2\pi\langle\xi,y\rangle}\hat u_0(\xi)$, where $\hat u_0\in C^\infty_c(\hat G)'$.

Mixed-Fourier-norm spaces of holomorphic functions are discussed in Section 6. Namely, we consider the space $\mathcal{A}_\mathcal{X}(G\times Y)$ of all holomorphic functions in the mixed-Fourier-norm space $\mathcal{X}(G\times Y)$,
$$
\mathcal{A}_\mathcal{X}(G\times Y)=\mathcal{X}(G\times Y)\cap\mathrm{Hol}(G\times Y).
$$
The principal result of this paper is Proposition 6, which shows that
$$
\mathfrak{F}:\mathcal{A}_\mathcal{X}(G\times Y)\to e^{-2\pi\langle\cdot,\cdot\rangle}\cdot\Xi(\hat G,\rho)
$$
is an isometry, where the Banach space $\Xi(\hat G,\rho)$ is a weighted analogue of $\Xi(\hat G)$. Under very mild assumptions, Fourier transform is in fact an isometric isomorphism between these spaces.

Section 7 is dedicated to the following questions. In practice, one is given a complex domain $\Omega$, with $G\times Y$ bi-holomorphically mapped into an open dense subset of $\Omega$. In the archetypical example of the unit disc $\Omega=\mathbb{D}$, the space $G\times Y=\mathbb{T}\times\mathbb{R}_+$ sits inside $\Omega$ in the form of the punctured disc $\mathbb{D}\setminus\{0\}$. In order to recover holomorphic functions on $\mathbb{D}$ from holomorphic functions on $\mathbb{D}\setminus\{0\}$ one needs boundedness in a neighbourhood of $0$. This boundedness is provided by Paley-Wiener-type results thanks to the right choice of the norm on the Fourier side. In the usual Paley-Wiener theorem on $\mathbb{R}$, Fourier transform maps an $L^2$ function with a holomorphic extension to the upper half-plane bounded around $+\imath\infty$ to an $L^2$ function supported in the positive semi-axis. Hidden in Proposition 6 above are similar properties, which are spelled out explicitly in this section. Namely, under mild conditions, we have
$$
\mathrm{supp}\,\hat u_0\subset\overline{\hat G_+},\quad\forall u\in\mathcal{A}_\mathcal{X}(G\times Y),
$$
where $\overline{\hat G_+}=\hat G\cap[0,+\infty)^n$ is the non-negative sector, and $\exists y_0>0$ such that
$$
\sup_{y>y_0}\|u(\cdot,y)\|_\infty<\infty,\quad\forall u\in\mathcal{A}_\mathcal{X}(G\times Y),
$$
where $y>y_0$ is understood as $y-y_0\in(\mathbb{R}_+)^n$.

In the final Section 8 we illustrate the theory on three notable examples: the elliptic model of the unit disc (as mentioned above), the parabolic model of the half-plane, and the hyperbolic model of the half-plane. These examples, although all coming from the geometric models of the unit disc, have remarkably different structure and features.

\section{Fourier transform on Abelian principal bundles $G\times Y$ and $G$-tempered distributions}

The aim of this section is to introduce definitions and notations regarding the (``half-'')Fourier transform on a trivial Abelian principal bundle $G\times Y$.

\subsection*{Functions on $G$}

For a locally compact Abelian group $G$, its Pontrjagin dual group $\hat G$ is another locally compact Abelian group defined in terms of unitary characters, but in case of an Abelian Lie group it can be given more explicitly. Our setting is geometric in nature, and we will be mostly interested in connected Abelian Lie groups and their dual groups. The model cases are $\mathbb{R}$, $\mathbb{Z}$ and $\mathbb{T}=\mathbb{R}/\mathbb{Z}$, equipped with their standard measures (Lebesgue or counting). All we need to remember is that
$$
\hat{\mathbb{R}}=\mathbb{R},\quad\hat{\mathbb{T}}=\mathbb{Z},\quad\hat{\mathbb{Z}}=\mathbb{T}.
$$
If $G$ is an Abelian Lie group, it can be thought of as
$$
G=\prod_{i=1}^nG_i,\quad G_i\in\{\mathbb{R},\mathbb{T},\mathbb{Z}\},\quad i=1,\ldots,n,\quad n\in\mathbb{N}_0,
$$
with (possibly discrete) group coordinates $x=(x_1,\ldots,x_n)\in G$, $x_i\in G_i$, $i=1,
\ldots,n$. The Haar measure, formally denoted by $dx=dx_1\ldots dx_n$, is the product of the standard measures on each component. The Pontrjagin dual group $\hat G$ is another such Abelian Lie group,
$$
\hat G=\prod_{i=1}^n\hat G_i,
$$
with group coordinates $\xi=(\xi_1,\ldots,\xi_n)\in\hat G$, $\xi_i\in\hat G_i$, $i=1,\ldots,n$. There is a natural additive bi-character $\langle\cdot,\cdot\rangle:G\times\hat G\to\mathbb{R}$ given by
$$
\langle x,\xi\rangle=\sum_{i=1}^nx_i\xi_i,\quad\forall(x,\xi)\in G\times\hat G,
$$
as well as a unitary bi-character
$$
\chi(x,\xi)=e^{2\pi\imath\langle x,\xi\rangle},\quad\forall(x,\xi)\in G\times\hat G.
$$
It is easily checked that
\begin{equation}
\partial_{x_i}\chi(x,\xi)=2\pi\imath\,\xi_i\,\chi(x,\xi),\quad\partial_{\xi_j}\chi(x,\xi)=2\pi\imath\,x_j\,\chi(x,\xi),\quad\forall x\in G,\quad\forall\xi\in\hat G,\label{chiprop}
\end{equation}
for those $i$ such that $G_i$ is connected and for those $j$ such that $\hat G_j$ is connected. For a function $f\in L^1(G)$ its Fourier transform $\hat f\in C(\hat G)$ is defined by
$$
\hat f(\xi)=\mathfrak{F}f(\xi)=\int\limits_Gf(x)\chi(x,\xi)^*dx=\int\limits_Gf(x)e^{-2\pi\imath\langle x,\xi\rangle}dx,\quad\forall\xi\in\hat G,
$$
where integration is understood in terms of the chosen Haar measure. If $\hat f\in L^1(\hat G)$ as well then a Fourier inversion formula holds,
$$
f(x)=\mathfrak{F}^{-1}\hat f(x)=\int\limits_{\hat G}\hat f(\xi)\chi(x,\xi)d\xi=\int\limits_{\hat G}\hat f(\xi)e^{2\pi\imath\langle x,\xi\rangle}d\xi,\quad\forall x\in G.
$$
Note that if we adopt the point of view of the group $\hat G$, then the Fourier transform of the integrable function $\hat f\in L^1(\hat G)$ will be
$$
\mathfrak{F}\hat f(x)=\int\limits_{\hat G}\hat f(\xi)\chi(\xi,x)^*d\xi=\int\limits_{\hat G}\hat f(\xi)e^{-2\pi\imath\langle\xi,x\rangle}d\xi=\breve f(x),\quad\forall x\in G,
$$
where $\breve f\in C(G)$ is given by $\breve f(x)=f(-x)$ for all $x\in G$, and $-x=(-x_1,\ldots,-x_n)$.

The Schwartz-Bruhat space $\mathcal{S}(G)$ of rapidly decaying test functions is a nuclear Fr\'echet space defined in terms of seminorms involving differential operators with polynomial coefficients. However, using nuclearity, we can construct them more explicitly by: setting $\mathcal{S}(\mathbb{R})$ to be the usual Schwartz space; letting $\mathcal{S}(\mathbb{T})=C^\infty(\mathbb{T})$; and
$$
\mathcal{S}(\mathbb{Z})=\left\{f:\mathbb{Z}\to\mathbb{C}\,\vline\quad\|f\|_m=\sup_{n\in\mathbb{Z}}|n^mf(n)|<\infty,\quad\forall m\in\mathbb{N}_0\right\},
$$
the space of rapidly decaying sequences; and defining
$$
\mathcal{S}(G)=\bigotimes_{i=1}^n\mathcal{S}(G_i).
$$
Tempered distributions $u\in\mathcal{S}(G)'$ are then continuous linear functionals.

It is clear that $\mathcal{S}(G)\subset L^1(G)$, and the Fourier transform $\mathfrak{F}:\mathcal{S}(G)\to\mathcal{S}(\hat G)$ is an isomorphism (see \cite{Osb75} and references therein). The Fourier transform $\hat u\in\mathcal{S}(\hat G)'$ of a distribution $u\in\mathcal{S}(G)'$ is defined by duality,
$$
\hat u(\hat f)=u(\breve f),\quad\forall f\in\mathcal{S}(G),
$$
since $f\mapsto\breve f$ is an automorphism of $\mathcal{S}(G)$. This is consistent, because if $u\in L^1(G)\cap L^2(G)$ then
$$
\hat u(\hat f)=\int\limits_{\hat G}\hat u(\xi)\hat f(\xi)d\xi=\int\limits_{G}u(x)\breve f(x)dx=u(\breve f),\quad\forall f\in\mathcal{S}(G),
$$
where we used $\mathfrak{F}\bar{\breve f}=\overline{\mathfrak{F}f}$ and Plancher\'el theorem. It is further shown that $\mathfrak{F}:\mathcal{S}(G)'\to\mathcal{S}(\hat G)'$ is again an isomorphism.

\subsection*{Functions on $G\times Y$}

Consider now a trivial principal $G$-bundle $G\times Y$ where $Y$ is a connected manifold (without boundary). The unique tensor product $\mathcal{D}(G\times Y)\doteq\mathcal{S}(G)\,\hat\otimes\,C^\infty_c(Y)$, a complete nuclear space, will be our space of test functions that are of compact support in $Y$ and of rapid decay in $G$. Since the tensor product $\mathcal{S}(G)\,\hat\otimes\,C^\infty_c(Y)$ of two nuclear spaces is inductive, we have that
\begin{equation}
\mathcal{S}(G)\,\hat\otimes\, C^\infty_c(Y)=\mathcal{S}(G)\,\hat\otimes\,\lim_{K\Subset Y}C^\infty_c(K)=\lim_{K\Subset Y}\mathcal{S}(G)\,\hat\otimes\, C^\infty_c(K).\label{IndLim}
\end{equation}
The dual space will consist of distributions $u\in\mathcal{D}(G\times Y)'=\mathcal{S}(G)'\,\hat\otimes\,C^\infty_c(Y)'$ which we will call $G$-tempered. We will define the Fourier transform $\mathfrak{F}:\mathcal{D}(G\times Y)\to\mathcal{D}(\hat G\times Y)$ by continuous extension of the componentwise operator $\mathfrak{F}\otimes\id:\mathcal{S}(G)\,\otimes\,C^\infty_c(Y)\to\mathcal{S}(\hat G)\,\otimes\,C^\infty_c(Y)$ acting on the vector space tensor product. It is clear that $\mathfrak{F}:\mathcal{D}(G\times Y)\to\mathcal{D}(\hat G\times Y)$ is an isomorphism, given explicitly by
$$
\hat f(\xi,y)=\int\limits_Gf(x,y)\chi(x,\xi)^*dx=\int\limits_Gf(x,y)e^{-2\pi\imath\langle x,\xi\rangle}dx,\quad\forall f\in\mathcal{D}(G\times Y),\quad\forall(\xi,y)\in\hat G\times Y.
$$
Note that due to the decay properties of $f$, the integral above converges absolutely, uniformly in $(\xi,y)\in K$ for every compact $K\Subset\hat G\times Y$. In particular, limits and differentiation can be performed unrestrictedly with respect to the variables $(\xi,y)$.

We define the Fourier transform of a $G$-tempered distribution by
$$
\hat u(\hat f)=u(\breve f),\quad\breve f(x,y)=f(-x,y),\quad\forall(x,y)\in G\times Y,\quad\forall f\in\mathcal{D}(G\times Y),\quad\forall u\in\mathcal{D}(G\times Y)'.
$$
This defines an isomorphism $\mathfrak{F}:\mathcal{D}(G\times Y)'\to\mathcal{D}(\hat G\times Y)'$.

One of the most important properties of Fourier transform is its behaviour with respect to certain differential operators. The classes $\mathrm{D}_G(G\times Y)$ and $\mathrm{D}_M(G\times Y)$ of differential operators appropriate for our analysis will be introduced in Section \ref{MultDiffOp}. At this moment we will make reference only to the classes $\mathrm{D}(Y)$ and $\mathrm{D}(G\times Y)$ of all differential operators on $Y$ and $G\times Y$, respectively, with smooth coefficients. The natural embedding $\mathrm{D}(Y)\subset\mathrm{D}(G\times Y)$ will be assumed by default.

Using properties (\ref{chiprop}) and integration by parts one can show that
\begin{eqnarray}
-2\pi\imath\,\widehat{[x_jf]}(\xi,y)=\partial_{\xi_j}\hat f(\xi,y),\quad\widehat{[\partial_{x_i}f]}(\xi,y)=2\pi\imath\,\xi_i\hat f(\xi,y),\nonumber\\
\widehat{[P(\partial_y)f]}(\xi,y)=P(\partial_y)\hat f(\xi,y),\quad\forall(\xi,y)\in\hat G\times Y,\quad\forall f\in\mathcal{D}(G\times Y),\label{^fprop}
\end{eqnarray}
for all $i$ such that $G_i$ is connected, all $j$ such that $\hat G_j$ is connected, and every differential operator $P(\partial_y)\in\mathrm{D}(Y)\subset\mathrm{D}(G\times Y)$. 
\section{Multipliers and differential operators acting on the test function space $\mathcal{D}(G\times Y)$}
\label{MultDiffOp}

The algebra of multipliers for the classical Schwartz space $\mathcal{S}(\mathbb{R})$ was found by L. Schwartz himself \cite{Sch78}, and is known to be
$$
\mathcal{O}_M(\mathbb{R})=\left\{u\in C^\infty(\mathbb{R})\,\vline\quad\left(\forall \alpha\in\mathbb{N}_0\right)\left(\exists m_\alpha^u\in\mathbb{N}_0\right)\quad\sup_{x\in\mathbb{R}}\left|(1+x^2)^{-m_\alpha^u}D_x^\alpha f(x)\right|<\infty\right\},
$$
equipped with the seminorms
$$
\|u\|_{\alpha,f}=\sup_{x\in\mathbb{R}}\left|f(x)D_x^\alpha u(x)\right|,\quad\forall u\in\mathcal{O}_M(\mathbb{R}),\quad\forall\alpha\in\mathbb{N}_0,\quad\forall f\in\mathcal{S}(\mathbb{R}).
$$
It was shown that the pointwise multiplication map
$$
\mathcal{O}_M(\mathbb{R})\times\mathcal{S}(\mathbb{R})\ni(u,f)\mapsto u\cdot f\in\mathcal{S}(\mathbb{R})
$$
is continuous in each variable separately (and even hypocontinuous). However, it was discovered that the multiplication map is not jointly continuous (see \cite{Lar13} for a nice review of such problems). Therefore the corresponding result for $\mathcal{S}(\mathbb{R}^2)$ cannot be reduced to that for $\mathcal{S}(\mathbb{R})$ by means of topological tensor product $\mathcal{O}_M(\mathbb{R}^2)=\mathcal{O}_M(\mathbb{R})\,\hat\otimes\,\mathcal{O}_M(\mathbb{R})$, and has to be proven from scratch. More so for our rather specific space $G\times Y$, which has different properties in different ``directions''. Therefore, for the sake of completeness, we will define our spaces and establish their basic properties below.

Denote by $\mathrm{D}_G(G\times Y)$ the algebra of differential operators $P(\partial)$ with smooth coefficients, acting on functions on $G\times Y$, which are $G$-invariant, i.e.,
$$
P(\partial_{x,y})[f(x+x_0,y)]=[P(\partial_{x,y})f](x+x_0,y),\quad\forall x,x_0\in G,\quad\forall y\in Y,\quad\forall f\in C^\infty(G\times Y).
$$
This is nothing else but an invariant description of differential operators with coefficients not depending on $x\in G$. The following representation will be used on several occasions in the sequel.

\begin{lemma}\label{PFactorLemma} A higher order Leibniz rule holds as follows:
$$
\left(\,\forall P(\partial)\in\mathrm{D}_G(G\times Y)\,\right)\quad\left(\,\exists N_P\in\mathbb{N}\,\right)\quad\left(\,\exists\{P_k(\partial),Q_k(\partial)\}_{k=1}^{N_P}\subset\mathrm{D}_G(G\times Y)^2\,\right)
$$
$$
P(\partial)[fh]=\sum_{k=1}^{N_P}[P_k(\partial)f]\cdot[Q_k(\partial)h],\quad\forall f,h\in C^\infty(G\times Y).
$$
\end{lemma}
\begin{proof} Let us first assert the above statement for the case where $Y\subset\mathbb{R}^m$ is an open subset. Then in coordinates $(x_1,\ldots,x_n,y_1,\ldots,y_m)\in G\times Y$ the differential operator $P(\partial)$ can be written as
$$
P(\partial)=\sum_{|\alpha|\le\operatorname{ord}P}a_\alpha(y)D^\alpha,
$$
where $\operatorname{ord}P$ is the order of $P(\partial)$, and the non-dependence of $a_\alpha$ on $x$ reflects the $G$-invariance of $P(\partial)$. A repeated application of Leibniz rule gives
$$
P(\partial)[fh]=\sum_{|\alpha|\le\operatorname{ord}P}a_\alpha(y)D^\alpha[fh]=\sum_{|\alpha|\le\operatorname{ord}P}a_\alpha(y)\sum_{|\beta|\le\operatorname{ord}P}{\alpha\choose\beta}[D^\beta f]\cdot[D^{\alpha-\beta}h]
$$
$$
=\sum_{|\beta|\le\operatorname{ord}P}[D^\beta f]\sum_{|\alpha|\le\operatorname{ord}P}a_\alpha(y){\alpha\choose\beta}\cdot[D^{\alpha-\beta}h]=\sum_{|\beta|\le\operatorname{ord}P}[P_\beta(\partial)f]\cdot[Q_\beta(\partial)h],
$$
where
$$
P_\beta(\partial)\doteq D^\beta,\quad Q_\beta(\partial)\doteq\sum_{|\alpha|\le\operatorname{ord}P}a_\alpha(y){\alpha\choose\beta}\cdot D^{\alpha-\beta},
$$
with the convention that
$$
{\alpha\choose\beta}=0,\quad\forall\beta\not\le\alpha.
$$
It remains to relabel $\beta\to k$ by enumeration with $N_P\doteq\#\{\beta\,|\,|\beta|<\operatorname{ord}P\}$ and to note that $P_k(\partial)$ and $Q_k(\partial)$ thus defined are $G$-invariant for $k=1,\ldots,N_P$.

Let us now consider a connected manifold $Y$ with $\dim Y=m\in\mathbb{N}$ and $P(\partial)\in\mathrm{D}_G(G\times Y)$. It is known that a connected smooth manifold has a finite atlas (see \cite{Sol74}, but a simpler argument can be made using triangulations). Let $\{U_j,\phi^j\}_{j=1}^N$ be a finite atlas for $Y$ with its subordinate smooth partition of unity $\{\chi_j\}_{j=1}^N$. Then for every fixed $j=1,\ldots,N$,
$$
\id\times\phi_j:G\times U_j\to G\times\phi_j(U_j)\subset G\times\mathbb{R}^m
$$
is a $G$-equivariant diffeomorphism and
$$
\tilde P^j(\partial)\doteq[\id\times\phi_j]_*^{-1}P(\partial)[\id\times\phi_j]_*\in\mathrm{D}_G(G\times\phi_j(U_j)),
$$
where $[\id\times\phi_j]_*:C^\infty(G\times\phi_j(U_j))\to C^\infty(G\times U_j)$ is the pullback map,
$$
\left([\id\times\phi_j]_*f\right)(x,y)=f(x,\phi_j(y)),\quad\forall(x,y)\in G\times U_j,\quad\forall f\in C^\infty(G\times\phi_j(U_j)).
$$
Thus, by the above statement,
$$
\tilde P^j(\partial)[fh]=\sum_{k=1}^{N_P}[\tilde P^j_k(\partial)f]\cdot[\tilde Q^j_k(\partial)h],\quad\forall f,h\in C^\infty(G\times\phi_j(U_j)),
$$
for some $\tilde P^j_k(\partial),\tilde Q^j_k(\partial)\in\mathrm{D}_G(G\times\phi_j(U_j))$, $k=1,\ldots, N_P$. Denote by $\eta_j\in C^\infty(G\times Y)$ the $G$-invariant function
$$
\eta_j(x,y)\doteq\frac{\chi_j(y)}{\sqrt{\sum\limits_{\ell=1}^N\chi_\ell(y)^2}},\quad\forall(x,y)\in G\times Y,
$$
and let
$$
P^j_k(\partial)\doteq\eta_j\cdot[\id\times\phi_j]_*\tilde P^j_k(\partial)[\id\times\phi_j]_*^{-1}\in\mathrm{D}_G(G\times Y),
$$
$$
Q^j_k(\partial)\doteq\eta_j\cdot[\id\times\phi_j]_*\tilde Q^j_k(\partial)[\id\times\phi_j]_*^{-1}\in\mathrm{D}_G(G\times Y),
$$
where the smooth extension beyond $G\times U_j$ by zero is understood thanks to $\supp\eta_j\subset U_j$. Set $f_j\doteq f|_{G\times U_j}\in C^\infty(G\times U_j)$ for every $f\in C^\infty(G\times Y)$ for convenience, and note that
$$
\left(\eta_j^2P(\partial)[fh]\right)|_{G\times U_j}=\eta_j^2P(\partial)[f_jh_j]=\eta_j^2[\id\times\phi_j]_*\tilde P^j(\partial)\left([\id\times\phi_j]_*^{-1}f_j\cdot[\id\times\phi_j]_*^{-1}h_j\right)
$$
$$
=\eta_j^2[\id\times\phi_j]_*\sum_{k=1}^{N_P}[\tilde P^j_k(\partial)[\id\times\phi_j]_*^{-1}f_j]\cdot[\tilde Q^j_k(\partial)[\id\times\phi_j]_*^{-1}h_j]
$$
$$
=\sum_{k=1}^{N_P}[P^j_k(\partial)f_j]\cdot[Q^j_k(\partial)h_j]=\left(\sum_{k=1}^{N_P}[P^j_k(\partial)f]\cdot[ Q^j_k(\partial)h]\right)\Bigr|_{G\times U_j},\quad\forall f,h\in C^\infty(G\times Y),
$$
while
$$
\left(\eta_j^2P(\partial)[fh]\right)|_{G\times U_j^\complement}=0=\left(\sum_{k=1}^{N_P}[P^j_k(\partial)f]\cdot[ Q^j_k(\partial)h]\right)\Bigr|_{G\times U_j^\complement},\quad\forall f,h\in C^\infty(G\times Y),
$$
which leads to the globally valid identity
$$
\eta_j^2P(\partial)[fh]=\sum_{k=1}^{N_P}[P^j_k(\partial)f]\cdot[Q^j_k(\partial)h],\quad\forall f,h\in C^\infty(G\times Y).
$$
Finally, after performing these constructions for all $j=1,\ldots,N$ and observing that $\{\eta_j^2\}_{j=1}^N$ is a partition of unity on $G\times Y$, we write
$$
P(\partial)[fh]=\sum_{j=1}^N\eta_j^2P(\partial)[fh]=\sum_{j=1}^N\sum_{k=1}^{N_P}[P^j_k(\partial)f]\cdot[Q^j_k(\partial)h],\quad\forall f,h\in C^\infty(G\times Y).
$$
It only remains to relabel $N\cdot N_P\mapsto N_P$ and $(j,k)\mapsto k$.
\end{proof}

Denote
\begin{equation}
\langle x\rangle\doteq1+\sum_{i\in\mathcal{I}_G}x_i^2,\quad\mathcal{I}_G\doteq\left\{i\in\{1,\ldots,n\}\,\vline\quad G_i\neq\mathbb{T}\right\},\quad\forall x\in G.\label{<.>DefEq}
\end{equation}
Define the seminorms
$$
\|f\|_{P,m,K}\doteq\sup_{(x,y)\in G\times K}\left|\langle x\rangle^mP(\partial)f(x,y)\right|,
$$
$$
\forall f\in C^\infty(G\times Y),\quad\forall P\in\mathrm{D}_G(G\times Y),\quad\forall m\in\mathbb{N}_0,\quad\forall K\Subset Y.
$$

\begin{lemma}\label{SeminormsLemma} For every compact $K\subset Y$, the topology of the space $\mathcal{S}(G)\hat\otimes\,C^\infty_c(K)$ can be given by the seminorms $\|.\|_{P,m,K}$
\end{lemma}
\begin{proof} It is sufficient to show that the original product topology restricted to the algebraic tensor product $\mathcal{S}(G_1)\otimes\ldots\otimes\mathcal{S}(G_n)\otimes C^\infty_c(K)$ is equivalent to the one given by seminorms $\|.\|_{P,m,K}$. It is known that the standard topology on $C^\infty_c(K)$ can be given by seminorms
$$
\|h\|_{P,K}\doteq\sup_{y\in K}\left|P(\partial)h\right|,\quad\forall h\in C^\infty(K),\quad\forall P(\partial)\in\mathrm{D}(Y),
$$
where $\mathrm{D}(Y)$ is the space of all differential operators with smooth coefficients acting on functions on $Y$. For every two multi-indices $\alpha,\beta\in\mathbb{N}_0^n$ such that $\alpha_i=0$ for $G_i=\mathbb{T}$ and $\beta_j=0$ for $G_j=\mathbb{Z}$,
$$
\sup_{x_1\in G_1}\left|x_1^{\alpha_1}\partial_{x_1}^{\beta_1}f_1(x_1)\right|\cdot\ldots\cdot\sup_{x_n\in G_n}\left|x_n^{\alpha_n}\partial_{x_n}^{\beta_n}f_n(x_n)\right|\cdot\|h\|_{P,K}
$$
$$
=\sup_{(x,y)\in G\times K}\left|x^\alpha D^\beta_x P(\partial_y)f_1(x_1)\ldots f_n(x_n)h(y)\right|
$$
$$
\le\sup_{(x,y)\in G\times K}\left|\langle x\rangle^{|\alpha|}D^\beta_x P(\partial_y)f_1(x_1)\ldots f_n(x_n)h(y)\right|=\|f_1\ldots f_nh\|_{D^\beta_xP,|\alpha|,K}\,,
$$
$$
\forall f_i\in\mathcal{S}(G_i),\quad i=1,\ldots,n,\quad\forall h\in C^\infty_c(K),\quad\forall P(\partial_y)\in\mathrm{D}(Y),
$$
where $D^\beta_xP\in\mathrm{D}_G(G\times Y)$ is understood in the obvious way. Conversely, let $P(\partial)\in\mathrm{D}_G(G\times Y)$ and $m\in\mathbb{N}_0$. A repeated application of Lemma \ref{PFactorLemma} gives
$$
P(\partial)[f_1\ldots f_nh]=\sum_{\alpha\in\mathbb{A}}[P_{\alpha_1}(\partial)f_1]\cdot\ldots\cdot[P_{\alpha_n}(\partial)f_n]\cdot[P_{\alpha_{n+1}}h]
$$
for some finite set of multi-indices $\mathbb{A}$ and $P_{\alpha_i}(\partial)\in\mathrm{D}_G(G\times Y)$, $i=1,\ldots,n+1$, $\alpha\in\mathbb{A}$. Set $m_i=m$ for $G_i\neq\mathbb{T}$ and $m_i=0$ for $G_i=\mathbb{T}$, so that
$$
\langle x\rangle^m\le\prod_{i=1}^n(1+x_i^2)^{m_i},\quad\forall x\in G.
$$
Then
$$
\|f_1\ldots f_nh\|_{P,m,K}=\sup_{(x,y)\in G\times K}\left|\langle x\rangle^mP(\partial)f_1(x_1)\ldots f_n(x_n)h(y)\right|
$$
$$
\le\sum_{\alpha\in\mathbb{A}}\Biggl[\sup_{(x,y)\in G\times K}\left|(1+x_1^2)^{m_1}P_{\alpha_1}(\partial)f_1(x_1)\right|\times\ldots\times\sup_{(x,y)\in G\times K}\left|(1+x_n^2)^{m_n}P_{\alpha_n}(\partial)f_n(x_n)\right|
$$
\begin{equation}
\times\sup_{(x,y)\in G\times K}\left|P_{\alpha_{n+1}}(\partial)h(y)\right|\Biggr].\label{Estimate}
\end{equation}
Note that due to the $G$-invariance of $P_{\alpha_i}(\partial)$, we have
$$
P_{\alpha_i}(\partial)f_i(x_i)=\sum_{\beta_i=1}^{N_i}a^i_{\beta_i}(y)\partial_{x_i}^{\beta_i}f_i(x_i),
$$
so that
$$
\sup_{(x,y)\in G\times K}\left|(1+x_i^2)^{m_i}P_{\alpha_i}(\partial)f_i(x_i)\right|\le\sum_{\beta_i=1}^{N_i}\sup_{y\in K}\left|a^i_{\beta_i}(y)\right|\cdot\sup_{x_i\in G_i}\left|(1+x_i^2)^{m_i}\partial_{x_i}^{\beta_i}f_i(x_i)\right|
$$
\begin{equation}
\le\sum_{\beta_i=1}^{N_i}\sum_{q_i=1}^{m_i}\sup_{y\in K}\left|a^i_{\beta_i}(y)\right|{m_i\choose q_i}\cdot\sup_{x_i\in G_i}\left|x_i^{2q_i}\partial_{x_i}^{\beta_i}f_i(x_i)\right|.\label{FormG}
\end{equation}
On the other hand, by $G$-invariance
\begin{equation}
\sup_{(x,y)\in G\times K}\left|P_{\alpha_{n+1}}(\partial)h(y)\right|=\sup_{y\in K}\left|P_{\alpha_{n+1}}(\partial)|_{x=0}h(y)\right|.\label{FormY}
\end{equation}
Substituting (\ref{FormG}) and (\ref{FormY}) into (\ref{Estimate}) we will estimate $\|.\|_{P,m,K}$ from above in terms of the original seminorms of the space $\mathcal{S}(G_1)\otimes\ldots\otimes\mathcal{S}(G_n)\otimes C^\infty_c(K)$. This concludes the proof.
\end{proof}

\begin{corollary}\label{PContCorollary} Every $P(\partial)\in\mathrm{D}_G(G\times Y)$ acts continuously on $\mathcal{D}(G\times Y)$.
\end{corollary}
\begin{proof} Let $K\Subset Y$ and let $f\in\mathcal{S}(G)\hat\otimes\,C^\infty_c(K)$. Then
$$
\|P(\partial)f\|_{Q,m,K}=\sup_{(x,y)\in G\times K}\left|\langle x\rangle^mQ(\partial)P(\partial)f(x,y)\right|=\|f\|_{QP,m,K},\quad\forall Q\in\mathrm{D}_G(G\times Y),\quad\forall m\in\mathbb{N}_0,
$$
which by Lemma \ref{SeminormsLemma} shows that $P(\partial)$ acts continuously on $\mathcal{S}(G)\hat\otimes C^\infty_c(K)$. By (\ref{IndLim}), it follows that $P(\partial)$ acts continuously on the inductive limit $\mathcal{D}(G\times Y)$.
\end{proof}

We are now ready to proceed to the multipliers on $\mathcal{D}(G\times Y)$. Define the space of smooth functions with slow growth in $G$ as
\begin{align}
&\mathcal{O}_\mathrm{M}(G\times Y)\doteq &\nonumber\\
&\left\{u\in C^\infty(G\times Y)\,\vline\,\left(\forall P(\partial)\in\mathrm{D}_G(G\times Y)\right)\left(\exists m_P^u\in\mathbb{N}_0\right)\left(\forall K\Subset Y\right)\,\sup_{(x,y)\in G\times K}\left|\langle x\rangle^{-m_P^u}P(\partial)u(x,y)\right|<\infty\right\},&\nonumber
\end{align}
equipped with the seminorms
$$
\|u\|_{P,f}\doteq\sup_{(x,y)\in G\times Y}\left|f(x,y)P(\partial)u(x,y)\right|,
$$
$$
\forall u\in C^\infty(G\times Y),\quad\forall P(\partial)\in\mathrm{D}_G(G\times Y),\quad\forall f\in\mathcal{D}(G\times Y).
$$

\begin{proposition}\label{MultContProp} The pointwise multiplication map
$$
\mathcal{O}_\mathrm{M}(G\times Y)\times\mathcal{D}(G\times Y)\ni(u,f)\mapsto uf\in\mathcal{D}(G\times Y)
$$
is continuous in each variable separately.
\end{proposition}
\begin{proof} First fix $u\in\mathcal{O}_\mathrm{M}(G\times Y)$. Let $K\Subset Y$ and $f\in\mathcal{S}(G)\,\hat\otimes\,C^\infty_c(K)$. Take $\forall P(\partial)\in\mathrm{D}_G(G\times Y)$ and $\forall m\in\mathbb{N}_0$. Then by Lemma \ref{PFactorLemma} we have some $N_P\in\mathbb{N}$ and $\{P_k(\partial),Q_k(\partial)\}_{k=1}^{N_P}\subset\mathrm{D}_G(G\times Y)^2$ such that
\begin{equation}
P(\partial)[uf]=\sum_{k=1}^{N_P}[P_k(\partial)u]\cdot[Q_k(\partial)f],\label{PFactorEq1}
\end{equation}
hence
$$
\|uf\|_{P,m,K}=\sup_{(x,y)\in G\times K}\left|\langle x\rangle^mP(\partial)[u(x,y)f(x,y)]\right|
$$
$$
\le\sum_{k=1}^{N_P}\sup_{(x,y)\in G\times K}\left|\langle x\rangle^m[P_k(\partial)u(x,y)]\cdot[Q_k(\partial)f(x,y)]\right|
$$
$$
\le\sum_{k=1}^{N_P}\sup_{(x,y)\in G\times K}\left|\langle x\rangle^{-m_{P_k}^u}P_k(\partial)u(x,y)\right|\cdot\|f\|_{Q_k,m+m_{P_k}^u,K}\,,
$$
where in the last step we used the definition of $\mathcal{O}_\mathrm{M}(G\times Y)$ to find integers $m_{P_k}^u$. This shows that the multiplication by $u$ acts continuously on $\mathcal{S}(G)\hat\otimes C^\infty_c(K)$. By (\ref{IndLim}), it acts continuously on the inductive limit $\mathcal{D}(G\times Y)$.

Let us now fix $f\in\mathcal{D}(G\times Y)$. Take $K\Subset Y$ such that $f\in\mathcal{S}(G)\,\hat\otimes\,C^\infty_c(K)$. Then $uf\in\mathcal{S}(G)\,\hat\otimes\,C^\infty_c(K)$ for all $u\in\mathcal{O}_\mathrm{M}(G\times Y)$. Take $\forall P(\partial)\in\mathrm{D}_G(G\times Y)$ and $\forall m\in\mathbb{N}_0$. By the same decomposition (\ref{PFactorEq1}) we find that
$$
\|uf\|_{P,m,K}=\sup_{(x,y)\in G\times Y}\left|\langle x\rangle^mP(\partial)[u(x,y)f(x,y)]\right|\le
$$
$$
\sum_{k=1}^{N_P}\sup_{(x,y)\in G\times Y}\left|[\langle x\rangle^mQ_k(\partial)f(x,y)]\cdot P_k(\partial)u(x,y)\right|=\sum_{k=1}^{N_P}\|u\|_{P_k,f_k},
$$
where $f_k\in\mathcal{D}(G\times Y)$ are defined by
$$
f_k(x,y)=\langle x\rangle^mQ_k(\partial)f(x,y),\quad\forall(x,y)\in G\times Y.
$$
That indeed $f_k\in\mathcal{D}(G\times Y)$ follows from Corollary \ref{PContCorollary} and the fact that $G\times Y\ni(x,y)\to\langle x\rangle^m\in\mathbb{R}$ is a function from $\mathcal{O}_\mathrm{M}(G\times Y)$, which was already shown to act on $\mathcal{D}(G\times Y)$ above. Thus, we have proven that the multiplication $(u,f)\mapsto uf$ is continuous in the first variable, too.
\end{proof}

We can define the product of a distribution $u\in\mathcal{D}(G\times Y)'$ with a multiplier $v\in\mathcal{O}_\mathrm{M}(G\times Y)$ by letting $vu\in\mathcal{D}(G\times Y)'$ be defined by
$$
[vu](f)=u(vf),\quad\forall f\in\mathcal{D}(G\times Y).
$$

\begin{corollary} The multiplication map
$$
\mathcal{O}_\mathrm{M}(G\times Y)\times\mathcal{D}(G\times Y)'\ni(u,v)\mapsto vu\in\mathcal{D}(G\times Y)'
$$
is continuous in each variable separately.
\end{corollary}
\begin{proof} Follows immediately from Proposition \ref{MultContProp} and duality.
\end{proof}

Now we can define a natural class of differential operators acting on $\mathcal{D}(G\times Y)$ larger than $\mathrm{D}_G(G\times Y)$. Namely, we introduce the algebra $\mathrm{D}_M(G\times Y)$ of differential operators with smooth coefficients of slow growth in $G$ acting on functions on $G\times Y$ as
\begin{equation}
\mathrm{D}_M(G\times Y)\doteq\langle\mathrm{D}_G(G\times Y),\mathcal{O}_M(G\times Y)\rangle,\label{D_MDef}
\end{equation}
that is, the smallest subalgebra of $\mathrm{D}(G\times Y)$ (the algebra of all differential operators with smooth coefficients acting on functions on $G\times Y$) containing $\mathrm{D}_G(G\times Y)$ and multipliers $\mathcal{O}_M(G\times Y)$ (considered as differential operators of order zero).

\begin{proposition} Every $P(\partial)\in\mathrm{D}_M(G\times Y)$ acts continuously on $\mathcal{D}(G\times Y)$.
\end{proposition}
\begin{proof} Follows immediately from Corollory \ref{PContCorollary} and Proposition \ref{MultContProp}.
\end{proof}

\section{Function spaces on $G\times Y$}

In this section we will discuss $G$-invariant measures on $G\times Y$, measurable functions and their relations with distributions, as well as weak derivatives. In the present work all manifolds will be equipped with their standard (Lebesgue) measurable space structures. All measurable functions on a manifold will be assumed Lebesgue measurable with respect to the standard measurable space structure. On the manifold $G\times Y$ it is natural to consider $G$-invariant measures $\mu$, i.e., $d\mu(x+x_0,y)=d\mu(x,y)$ for all fixed $x_0\in G$.

\begin{lemma} For every $G$-invariant Radon measure $\mu$ on $G\times Y$ there exists a Radon measure $\nu$ on $Y$ such that $d\mu(x,y)=dx\,d\nu(y)$.
\end{lemma}
\begin{proof} For every $h\in C_c(Y)$ we define the invariant Radon measure $\mu_h$ on $G$ by $\mu_h(f)=\mu(hf)$ for all $f\in C_c(G)$. By the uniqueness of the Haar measure on $G$, there exists a positive linear functional $\nu$ on $C_c(Y)$ such that $\mu_h(f)=\mu(hf)=\nu(h)\int\limits_Gfdx$ for all $f\in C_c(G)$ and $h\in C_c(Y)$. The functional $\nu$ is also continuous, and hence a Radon measure on $Y$.
\end{proof}
The decomposition of the invariant measure $\mu$ allows us to equip $\hat G\times Y$ with an invariant measure $\hat\mu$ given by $d\hat\mu(\xi,y)=d\xi\,d\nu(y)$.

Given a $G$-invariant measure $\mu$ on $G\times Y$, we can speak of Lebesgue spaces $L^p(G\times Y)=L^p(G\times Y,\mu)$, $p\in[0,+\infty]$. The Hilbert space $L^2(G\times Y)$ is given by the inner product
$$
(f,g)_2=(f,g)_{2,\mu}=\int\limits_{G\times Y}\bar f(x,y)g(x,y)d\mu(x,y)=\int\limits_{G\times Y}\bar f(x,y)g(x,y)dx\,d\nu(y),
$$
for all $f,g\in L^0(G\times Y)$ for which the integral makes sense. For $p\in[1,+\infty]$ we will use $\|.\|_p=\|.\|_{p,\mu}$ to denote the usual norms.

Recall the definition of $\langle x\rangle$ from (\ref{<.>DefEq}) and let $d\mu_m(x,y)=\langle x\rangle^md\mu(x,y)$ define a sequence of (non-$G$-invariant) measures $\{\mu_m\}_{m=-\infty}^\infty$ on $G\times Y$. Every $u\in L^1(G\times Y,\mu_m)$, $m\in\mathbb{Z}$, gives rise to a distribution $u\in\mathcal{D}(G\times Y)'$ by
$$
u(f)=\int\limits_{G\times Y}u(x,y)f(x,y)d\mu(x,y),\quad\forall f\in\mathcal{D}(G\times Y).
$$
If $\supp\nu=Y$ (i.e., $\supp\mu_m=G\times Y$, $m\in\mathbb{Z}$) then this map $L^1(G\times Y,\mu_m)\to \mathcal{D}(G\times Y)'$ is injective. We have also
$$
|u(f)|\le\|u\|_{1,\mu_m}\cdot\|f\|_{1,\max\{0,m\},K}\,,\quad\forall f\in\mathcal{D}(G\times Y),\quad\forall K\Subset Y,\quad\supp f\subset K,
$$
which shows that the embedding
$$
L^1(G\times Y,\mu_m)\hookrightarrow\mathcal{D}(G\times Y)'
$$
thus defined is continuous. On the other hand,
$$
\mathcal{D}(G\times Y)\hookrightarrow L^p(G\times Y),\quad\forall p\in(0,+\infty],
$$
are continuous embeddings, too. Indeed, for $p<\infty$ we have
$$
\|f\|_p^p=\int\limits_{G\times K}d\mu(x,y)\frac1{\langle x\rangle^n}\cdot\langle x\rangle^n|f(x,y)|^p\le\nu(K)\int\limits_G\frac{dx}{\langle x\rangle^n}\cdot\left(\sup_{(x,y)\in G\times K}\left|\langle x\rangle^{\lceil\frac{n}{p}\rceil}f(x,y)\right|\right)^p
$$
$$
=\nu(K)\int\limits_G\frac{dx}{\langle x\rangle^n}\cdot\|f\|_{1,\lceil\frac{n}{p}\rceil,K}^p,\quad\forall f\in\mathcal{D}(G\times Y),\quad\forall K\Subset Y,\quad\supp f\subset K,
$$
while for $p=\infty$ simply
$$
\|f\|_\infty=\|f\|_{1,0,K},\quad\forall f\in\mathcal{D}(G\times Y),\quad\forall K\Subset Y,\quad\supp f\subset K.
$$

Let us now turn to the Fourier transform of measurable functions. For every measurable function $f:G\times Y\to\mathbb{C}$,
$$
\left(\forall y\in Y\right)\quad f(\cdot,y)\in L^1(G)\quad\Rightarrow\quad\hat f(\cdot,y)\in C(\hat G).
$$
Moreover, by Tonelli's and Plancher\'el's theorems,
$$
\|f\|_2^2=\int\limits_{G\times Y}|f(x,y)|^2d\mu(x,y)=\int\limits_Yd\nu(y)\int\limits_G|f(x,y)|^2dx=\int\limits_Yd\nu(y)\int\limits_{\hat G}|\hat f(\xi,y)|^2d\xi
$$
$$
=\int\limits_{\hat G\times Y}|\hat f(\xi,y)|^2d\hat\mu(\xi,y)=\|\hat f\|_2^2,\quad\forall f\in L^0(G\times Y),
$$
which is a variant of Plancher\'el theorem for the Fourier transform $\mathfrak{F}:L^2(G\times Y)\to L^2(\hat G\times Y)$. If $u\in L^1(G\times Y)\cap L^2(G\times Y)$ then the (strong) Fourier transform of $u$ as a function is consistent with its (weak) Fourier transform as a distribution,
$$
\hat u(\hat f)=u(\breve f)=\int\limits_{G\times Y}u(x,y)\breve f(x,y)d\mu(x,y)=\int\limits_{\hat G\times Y}\hat u(\xi,y)\hat f(\xi,y)d\hat\mu(x,y),\quad\forall f\in\mathcal{D}(G\times Y).
$$

Henceforth we will assume that the measure $\nu$ is given by a smooth, nowhere vanishing positive density. The smoothness and non-vanishing of $\nu$ (and thus, $\mu$) are imposed in order that the transpose operation on differential operators preserves the smoothness of coefficients. (They are also essential for the correspondence between the regularities of a locally integrable function and the corresponding distribution.) Namely, for every differential operator $P(\partial)\in\mathrm{D}(G\times Y)$ with smooth coefficients, acting on functions on $G\times Y$, there exists a unique differential operator $P(\partial)^\top\in\mathrm{D}(G\times Y)$ of the same kind, called $P(\partial)$-s transpose, such that
$$
\int\limits_{G\times Y}h(x,y)P(\partial)f(x,y)d\mu(x,y)=\int\limits_{G\times Y}f(x,y)P(\partial)^\top h(x,y)d\mu(x,y),\quad\forall f,h\in C^\infty_c(G\times Y).
$$
The map $P(\partial)\mapsto P(\partial)^\top$ is a linear involution of the space $\mathrm{D}(G\times Y)$, and
$$
\left[P(\partial)Q(\partial)\right]^\top=Q(\partial)^\top P(\partial)^\top,\quad\forall P(\partial),Q(\partial)\in\mathrm{D}(G\times Y).
$$
By $G$-invariance of $\mu$, we have
$$
P(\partial)\in\mathrm{D}_G(G\times Y)\quad\Leftrightarrow\quad P(\partial)^\top\in\mathrm{D}_G(G\times Y).
$$
On the other hand, if $Q(\partial)\in\mathrm{D}(G\times Y)$ is of order zero, such as a multiplier, then $Q(\partial)^\top=Q(\partial)$. Then it follows from the definition that
$$
P(\partial)\in\mathrm{D}_M(G\times Y)\quad\Leftrightarrow\quad P(\partial)^\top\in\mathrm{D}_M(G\times Y),
$$
that is, the operation of transposition is a linear involution on the space $\mathrm{D}_M(G\times Y)$. Note that every $P(\partial_y)\in\mathrm{D}(Y)\subset\mathrm{D}_G(G\times Y)$ has a transpose $P(\partial_y)^{\top_Y}\in\mathrm{D}(Y)\subset\mathrm{D}_G(G\times Y)$ with respect to $(Y,\nu)$, i.e.,
$$
\int\limits_Yh(y)P(\partial_y)f(y)d\nu(y)=\int\limits_Yf(y)P(\partial_y)^{\top_Y}h(y)d\nu(y),\quad\forall f,h\in C^\infty_c(Y).
$$
Then
$$
\int\limits_{G\times Y}h_1(x)h_2(y)P(\partial_y)f_1(x)f_2(y)dxd\nu(y)=\int\limits_{G\times Y}f_1(x)f_2(y)P(\partial_y)^{\top_Y}h_1(x)h_2(y)dxd\nu(y),
$$
$$
\forall f_1,h_1\in C^\infty_c(G),\quad\forall f_2,h_2\in C^\infty_c(Y),
$$
which by uniqueness shows that $P(\partial_y)^{\top_Y}=P(\partial_y)^\top$. Thus,
\begin{equation}
P(\partial)\in\mathrm{D}(Y)\quad\Leftrightarrow\quad P(\partial)^\top=P(\partial)^{\top_Y}\in\mathrm{D}(Y).\label{D_Y_TransposeEq}
\end{equation}

Weak derivatives of distributions are defined in terms of transposition. Namely, for a differential operator $P(\partial)\in\mathrm{D}_M(G\times Y)$, its action $P(\partial):\mathcal{D}(G\times Y)'\to\mathcal{D}(G\times Y)'$ is defined by
$$
[P(\partial)u](f)\doteq u\left(P(\partial)^\top f\right),\quad\forall f\in\mathcal{D}(G\times Y).
$$
Since $P(\partial)^\top:\mathcal{D}(G\times Y)\to\mathcal{D}(G\times Y)$ is continuous, so is $P(\partial):\mathcal{D}(G\times Y)'\to\mathcal{D}(G\times Y)'$ in the dual topology.

\begin{remark} Everything stated above for $G\times Y$ holds mutatis mutandis for $\hat G\times Y$.
\end{remark}

\section{Mixed-Fourier-norm spaces on $G\times Y$}\label{MixedFourierNormSection}

Now we will introduce the main function spaces of our study. Let $(\mathcal{Y}(Y),\|.\|_\mathcal{Y})$ be a continuously embedded Banach space consisting of distributions on $Y$, i.e.,
$$
\mathcal{Y}(Y)\hookrightarrow C^\infty_c(Y)'.
$$
By setting $\|v\|_\mathcal{Y}=+\infty$ if $\|v\|_\mathcal{Y}$ is undefined, we can assume without loss of generality that the norm $\|v\|_\mathcal{Y}$ is defined for all $v\in C^\infty_c(Y)'$. We will say that the embedding $\mathcal{Y}(Y)\hookrightarrow C^\infty_c(Y)'$ is uniform if
\begin{equation}
|v(h)|\le\wp(h)\cdot\|v\|_\mathcal{Y},\quad\forall h\in C^\infty_c(Y),\quad\forall v\in\mathcal{Y}(Y),\label{UniformEmbed}
\end{equation}
where $\wp\in C(C^\infty_c(Y),[0,+\infty))$, with $\wp(0)=0$, does not depend on $v$.

Let $(\Xi(\hat G),\|.\|_\Xi)$ be a continuously embedded Banach space consisting of locally integrable tempered distributions on $\hat G$, i.e.,
$$
\Xi(\hat G)\hookrightarrow\mathcal{S}(\hat G)',\quad\Xi(\hat G)\subset L^1_\mathrm{loc}(\hat G)\cap\mathcal{S}(\hat G)'.
$$
By setting $\|\hat u\|_\Xi=+\infty$ if $\|\hat u\|_\Xi$ is undefined, we can assume without loss of generality that the norm $\|\hat u\|_\Xi$ is defined for all $\hat u\in C^\infty_c(\hat G)'$. We will say that the embedding $\Xi(\hat G)\hookrightarrow\mathcal{S}(\hat G)'$ is uniform if
\begin{equation}
|\hat u(\hat f)|\le\varsigma(\hat f)\cdot\|\hat u\|_\Xi,\quad\forall\hat f\in\mathcal{S}(\hat G),\quad\forall\hat u\in\Xi(\hat G),\label{XiUniformEmbed}
\end{equation}
where $\varsigma\in C(\mathcal{S}(\hat G),[0,+\infty))$, with $\varsigma(0)=0$, does not depend on $\hat u$.

We will say that the norm $\|\cdot\|_\Xi$ has the quasi-lattice property if $\exists K_\Xi>0$ such that
\begin{equation}
\left(\forall\hat u_1,\hat u_2\in C^\infty_c(\hat G)'\right)\quad\hat u_1\ge0\quad\wedge\quad\hat u_2-\hat u_1\ge0\quad\Rightarrow\quad\|\hat u_1\|_\Xi\le K_\Xi\|\hat u_2\|_\Xi,\label{LatticeXiNorm}
\end{equation}
and that the Banach space $\Xi(\hat G)$ has the lattice property if
\begin{equation}
\left(\forall\hat u_1,\hat u_2\in C^\infty_c(\hat G)'\right)\quad\hat u_1\ge0\quad\wedge\quad\hat u_2-\hat u_1\ge0\quad\wedge\quad\hat u_2\in\Xi(\hat G)\quad\Rightarrow\quad\hat u_1\in\Xi(\hat G).\label{LatticeXiSpace}
\end{equation}

Consider the following exhaustion conditions on the Banach spaces $\Xi(\hat G)$ and $\mathcal{Y}(Y)$:
\begin{equation}
\Xi(\hat G)=\left\{\hat u\in C^\infty_c(\hat G)'\,\vline\quad\|\hat u\|_\Xi<\infty\right\},\label{XiNormExhaust}
\end{equation}
\begin{equation}
\mathcal{Y}(Y)=\left\{v\in C^\infty_c(Y)'\,\vline\quad\|v\|_\mathcal{Y}<\infty\right\}\label{YNormExhaust}.
\end{equation}

\begin{remark}\label{ExhaustionRemark} Note that if (\ref{XiNormExhaust}) is true then the quasi-lattice property (\ref{LatticeXiNorm}) implies the lattice property (\ref{LatticeXiSpace}). Moreover, in that case it is sufficient to verify the quasi-lattice property (\ref{LatticeXiNorm}) on locally integrable functions. The lattice property (\ref{LatticeXiSpace}) can be verified on locally integrable functions irrespective of (\ref{XiNormExhaust}).
\end{remark}
Indeed, either (\ref{LatticeXiNorm}) with (\ref{XiNormExhaust}), or (\ref{LatticeXiSpace}) alone, is non-trivial only if $\hat u_2\in\Xi(\hat G)\subset L^1_\mathrm{loc}(\hat G)$. The positive distributions $\hat u_1,\hat u_2$ are necessarily measures (e.g., Theorem 2.1.7 in \cite{HorI}), and $\hat u_2$ is absolutely continuous with respect to $d\xi$. Then from $0\le\hat u_1\le\hat u_2$ we find that $\hat u_1$ is also absolutely continuous with respect to $d\xi$, and is therefore locally integrable.

Let $L^1_\mathrm{loc}(\hat G,\mathcal{Y}(Y))$ be the space of locally Bochner-integrable functions on $\hat G$ taking values in the Banach space $\mathcal{Y}(Y)$. Then for every Bochner-measurable function $\hat u:\hat G\to\mathcal{Y}(Y)$,
\begin{equation}
\hat u\in L^1_\mathrm{loc}(\hat G,\mathcal{Y}(Y))\quad\Leftrightarrow\quad\|\hat u(\cdot)\|_\mathcal{Y}\in L^1_\mathrm{loc}(\hat G).\label{L1locEquiv}
\end{equation}
Denote
$$
\Xi(\hat G,\mathcal{Y}(Y))\doteq\left\{\hat u\in L^1_\mathrm{loc}(\hat G,\mathcal{Y}(Y))\,\vline\quad\|\hat u(\cdot)\|_\mathcal{Y}\in\Xi(\hat G)\right\}
$$
and
$$
\|\hat u\|_{\mathcal{Y},\Xi}\doteq\|\|\hat u(\cdot)\|_\mathcal{Y}\|_\Xi,\quad\forall\hat u\in L^1_\mathrm{loc}(\hat G,\mathcal{Y}(Y)).
$$
By definition, $\Xi(\hat G,\mathcal{Y}(Y))$ is a $\mathbb{C}$-conic space and $\|\cdot\|_{\mathcal{Y},\Xi}$ is a positive homogeneous function. The space $\Xi(\hat G,\mathcal{Y}(Y))$ is given the coarsest topology such that $\|\cdot-v\|_{\mathcal{Y},\Xi}$ is continuous for every $v\in\Xi(\hat G,\mathcal{Y}(Y))$. This makes $\Xi(\hat G,\mathcal{Y}(Y))$ a topological conic space.

\begin{lemma}\label{LatticeLemma} The following statements hold:
\begin{itemize}

\item[1.] If the lattice property (\ref{LatticeXiSpace}) is true then $\Xi(\hat G,\mathcal{Y}(Y))\subset L^1_\mathrm{loc}(\hat G,\mathcal{Y}(Y))$ is a vector subspace.

\item[2.] If the quasi-lattice property (\ref{LatticeXiNorm}) is true then $\|\cdot\|_{\mathcal{Y},\Xi}$ is a quasi-norm with the same triangle constant $K_\Xi$.

\end{itemize}
\end{lemma}
\begin{proof} Let $\hat u_1,\hat u_2\in\Xi(\hat G,\mathcal{Y}(Y))$ and $\lambda\in\mathbb{C}$. Then $\lambda\hat u_1+\hat u_2\in L^1_\mathrm{loc}(\hat G,\mathcal{Y}(Y))$ and
$$
0\le\|[\lambda\hat u_1+\hat u_2](\cdot)\|_\mathcal{Y}\le|\lambda|\|\hat u_1(\cdot)\|_\mathcal{Y}+\|\hat u_2(\cdot)\|_\mathcal{Y}\in\Xi(\hat G),
$$
which in view of (\ref{LatticeXiSpace}) implies $\|[\lambda\hat u_1+\hat u_2](\cdot)\|_\mathcal{Y}\in\Xi(\hat G)$ and thus, $\lambda\hat u_1+\hat u_2\in\Xi(\hat G,\mathcal{Y}(Y))$. This proves statement 1. To prove statement 2., the only non-trivial part is the triangle inequality. Let $\hat u_1,\hat u_2\in L^1_\mathrm{loc}(\hat G,\mathcal{Y}(Y))$. Then
$$
0\le\|[\hat u_1+\hat u_2](\cdot)\|_\mathcal{Y}\le\|\hat u_1(\cdot)\|_\mathcal{Y}+\|\hat u_2(\cdot)\|_\mathcal{Y},
$$
which thanks to (\ref{LatticeXiNorm}) gives
$$
\|\hat u_1+\hat u_2\|_{\mathcal{Y},\Xi}=\|\|[\hat u_1+\hat u_2](\cdot)\|_\mathcal{Y}\|_\Xi\le K_\Xi\|\|\hat u_1(\cdot)\|_\mathcal{Y}+\|\hat u_2(\cdot)\|_\mathcal{Y}\|_\Xi
$$
$$
\le K_\Xi\left(\|\|\hat u_1(\cdot)\|_\mathcal{Y}\|_\Xi+\|\|\hat u_2(\cdot)\|_\mathcal{Y}\|_\Xi\right)=K_\Xi\left(\|\hat u_1\|_{\mathcal{Y},\Xi}+\|\hat u_1\|_{\mathcal{Y},\Xi}\right),
$$
which completes the proof.
\end{proof}

Setting $\|\hat u\|_{\mathcal{Y},\Xi}=+\infty$ for $\hat u\in C^\infty_c(\hat G\times Y)'\setminus L^1_\mathrm{loc}(\hat G,\mathcal{Y}(Y))$ we can assume without loss of generality that $\|\hat u\|_{\mathcal{Y},\Xi}$ is defined for all $\hat u\in C^\infty_c(\hat G\times Y)'$.

\begin{lemma}\label{EmbeddingLemma} Suppose that the uniform embedding (\ref{UniformEmbed}) is true. Then

\begin{itemize}

\item[1.] A continuous embedding
\begin{equation}
L^1_\mathrm{loc}(\hat G,\mathcal{Y}(Y))\hookrightarrow C^\infty_c(\hat G\times Y)'\label{L1locEmbed}
\end{equation}
is given by
$$
\hat u(\hat f)=\int\limits_{\hat G}\hat u(\xi)\left(\hat f(\xi,.)\right)d\xi,\quad\forall\hat u\in L^1_\mathrm{loc}(\hat G,\mathcal{Y}(Y)),\quad\forall\hat f\in C^\infty_c(\hat G\times Y).
$$

\item[2.] The above injection restricts to a continuous embedding of topological conic spaces
$$
\Xi(\hat G,\mathcal{Y}(Y))\hookrightarrow\mathcal{D}(\hat G\times Y)'.
$$

\end{itemize}
\end{lemma}
\begin{proof} For the statement 1. note that from
$$
|\hat u(\hat fh)|\le\int\limits_{\hat G}|\hat u(\xi)(h)||\hat f(\xi)|d\xi\le\wp(h)\int\limits_{\hat G}\|\hat u(\xi)\|_\mathcal{Y}|\hat f(\xi)|d\xi
$$
$$
\le\wp(h)\|\hat f\|_\infty\int\limits_{K}\|\hat u(\xi)\|_\mathcal{Y}d\xi,\quad\forall K\Subset\hat G,\quad\forall\hat f\in C^\infty_c(\hat G),\quad\supp\hat f\subset K,\quad\forall h\in C^\infty_c(Y),
$$
it follows that the embedding (\ref{L1locEmbed}) is indeed continuous. For the statement 2. observe that for $\forall\hat u\in L^1_\mathrm{loc}(\hat G,\mathcal{Y}(Y))$, $\forall\hat f\in\mathcal{S}(\hat G)$ and $\forall h\in C^\infty_c(Y)$,
$$
|\hat u(\hat fh)|\le\int\limits_{\hat G}|\hat f(\xi)||\hat u(\xi)(h)|d\xi\le\wp(h)\int\limits_{\hat G}|\hat f(\xi)|\|\hat u(\xi)\|_\mathcal{Y}d\xi\le\wp(h)\int\limits_{\hat G}\hat g(\xi)\|\hat u(\xi)\|_\mathcal{Y}d\xi
$$
$$
=\wp(h)\bigl|\|\hat u(\cdot)\|_\mathcal{Y}(\hat g)\bigr|\le\wp(h)\|\|\hat u(\cdot)\|_\mathcal{Y}\|_\Xi\cdot C_{\hat g},
$$
where $\hat g\in\mathcal{S}(\hat G)$ is such that $|\hat f|\le g$, and $C_{\hat g}$ is a positive number given by the continuous embedding $\Xi(\hat G)\hookrightarrow\mathcal{S}(\hat G)'$.
\end{proof}

Consider the following conditions on the Banach space $\Xi(\hat G)$ and the norm $\|\cdot\|_\Xi$:
\begin{equation}
\||\hat u|\|_\Xi\le C_\Xi\,\|\hat u\|_\Xi,\quad\forall\hat u\in L^1_\mathrm{loc}(\hat G),\quad C_\Xi>0,\label{Xi|u|<u}
\end{equation}
\begin{equation}
\|\hat u\|_\Xi\le c_\Xi\,\||\hat u|\|_\Xi,\quad\forall\hat u\in L^1_\mathrm{loc}(\hat G),\quad c_\Xi>0,\label{Xiu<|u|}
\end{equation}
\begin{equation}
\hat u\in\Xi(\hat G)\quad\Rightarrow\quad|\hat u|\in\Xi(\hat G),\quad\forall\hat u\in L^1_\mathrm{loc}(\hat G),\label{u_|u|}
\end{equation}
\begin{equation}
|\hat u|\in\Xi(\hat G)\quad\Rightarrow\quad\hat u\in\Xi(\hat G),\quad\forall\hat u\in L^1_\mathrm{loc}(\hat G).\label{|u|_u}
\end{equation}

\begin{remark} Note that (\ref{Xi|u|<u}) combined with the exhaustion property (\ref{XiNormExhaust}) implies (\ref{u_|u|}), while (\ref{Xiu<|u|}) combined with the exhaustion property (\ref{XiNormExhaust}) implies (\ref{|u|_u}).
\end{remark}

\begin{lemma}\label{YXiNormEstLemma} Suppose that (\ref{XiNormExhaust}) and (\ref{YNormExhaust}) are true. Then (\ref{Xi|u|<u}) implies
\begin{equation}
\|\hat u\cdot v\|_{\mathcal{Y},\Xi}\le C_\Xi\,\|\hat u\|_\Xi\,\|v\|_\mathcal{Y},\quad\forall\hat u\in C^\infty_c(\hat G)',\quad\forall v\in C^\infty_c(Y)',\label{YXi<Xi*Y}
\end{equation}
while (\ref{Xiu<|u|}) implies
\begin{equation}
\|\hat u\cdot v\|_{\mathcal{Y},\Xi}\ge c_\Xi\,\|\hat u\|_\Xi\,\|v\|_\mathcal{Y},\quad\forall\hat u\in C^\infty_c(\hat G)',\quad\forall v\in C^\infty_c(Y)'.\label{YXi>Xi*Y}
\end{equation}
\end{lemma}
\begin{proof} Let (\ref{XiNormExhaust}), (\ref{YNormExhaust}) and (\ref{Xi|u|<u}) hold. Then (\ref{YXi<Xi*Y}) needs a verification only when $0<\|\hat u\|_\Xi\,\|v\|_\mathcal{Y}<\infty$, i.e., $0\neq\hat u\in\Xi(\hat G)$ and $0\neq v\in\mathcal{Y}(Y)$. In that case
$$
\|\hat u\cdot v\|_{\mathcal{Y},\Xi}=\|\|[\hat u\cdot v](\cdot)\|_\mathcal{Y}\|_\Xi=\||\hat u|\cdot\|v\|_\mathcal{Y}\|_\Xi=\||\hat u|\|_\Xi\,\|v\|_\mathcal{Y}\le C_\Xi\,\|\hat u\|_\Xi\,\|v\|_\mathcal{Y}.
$$
Now let (\ref{XiNormExhaust}), (\ref{YNormExhaust}) and (\ref{Xiu<|u|}) hold. Then (\ref{YXi>Xi*Y}) needs a verification only when $0<\|\hat u\cdot v\|_{\mathcal{Y},\Xi}<\infty$, i.e., $0\neq\hat u\cdot v\in L^1_\mathrm{loc}(\hat G,\mathcal{Y}(Y))$ and $\|[\hat u\cdot v](\cdot)\|_\mathcal{Y}\in\Xi(\hat G)$. This implies $v\in\mathcal{Y}(Y)$, $\hat u\in L^1_\mathrm{loc}(\hat G)$ and $|\hat u|\in\Xi(\hat G)$, whence $\|\hat u\|_\Xi\le c_\Xi\,\||\hat u|\|_\Xi<\infty$ and therefore $\hat u\in\Xi(\hat G)$. Now
$$
\|\hat u\cdot v\|_{\mathcal{Y},\Xi}=\|\|[\hat u\cdot v](\cdot)\|_\mathcal{Y}\|_\Xi=\||\hat u|\|_\Xi\,\|v\|_\mathcal{Y}\ge c_\Xi\,\|\hat u\|_\Xi\,\|v\|_\mathcal{Y},
$$
as desired.
\end{proof}

\begin{corollary}\label{XiEmbdedCorr} If (\ref{XiNormExhaust}), (\ref{YNormExhaust}) and (\ref{Xiu<|u|}) are true, then the natural injection $C^\infty_c(\hat G)'\,\otimes\, C^\infty_c(Y)'\hookrightarrow C^\infty_c(\hat G\times Y)'$ restricts to a continuous embedding of topological conic spaces
$$
\Xi(\hat G)\otimes\mathcal{Y}(Y)\hookrightarrow\Xi(\hat G,\mathcal{Y}(Y)).
$$
\end{corollary}
\begin{proof} Follows immediately from Lemma \ref{YXiNormEstLemma}.
\end{proof}

\begin{remark} Even when the continuous embedding in the above corollary is valid, it is still unknown if its image is dense.
\end{remark}

\begin{remark} If $\Xi(\hat G,\mathcal{Y}(Y))\subset L^1_\mathrm{loc}(\hat G,\mathcal{Y}(Y))$ is a vector subspace then the embeddings in the statement 2. of Lemma \ref{EmbeddingLemma} and in Corollary \ref{XiEmbdedCorr} are linear.
\end{remark}

If $\|\cdot\|_{\mathcal{Y},\Xi}$ is a quasi-norm, then the entourages
$$
\left\{(\hat u,\hat v)\in\Xi(\hat G,\mathcal{Y}(Y))^2\,\vline\quad\|\hat u-\hat v\|_{\mathcal{Y},\Xi}<\epsilon\right\}
$$
form a base of uniformity for $\Xi(\hat G,\mathcal{Y}(Y))$, making it a uniform conic space. In that case the embeddings in the statement 2. of Lemma \ref{EmbeddingLemma} and in Corollary \ref{XiEmbdedCorr} are uniformly continuous. Whether or not the uniform space $\Xi(\hat G,\mathcal{Y}(Y))$ is complete is a question that we will not be able to settle here. Instead, we will prove below a weaker partial result in that direction.

\begin{proposition}\label{XiMeasureClosed} Suppose that the quasi-lattice property (\ref{LatticeXiNorm}) holds, so that, by the statement 2. of Lemma \ref{LatticeLemma}, $\Xi(\hat G,\mathcal{Y}(Y))$ is a uniform conic space. Suppose that the embedding
\begin{equation}
\Xi(\hat G)\hookrightarrow L^0(\hat G)\label{XiMeasureProp}
\end{equation}
is continuous. Assume further that (\ref{Xiu<|u|}) is true. If $\{\hat u_m\}_{m=1}^\infty\subset\Xi(\hat G,\mathcal{Y}(Y))$ is a Cauchy sequence, then $\exists\hat u\in\Xi(\hat G,\mathcal{Y}(Y))$ such that $\hat u_m\xrightarrow[m\to\infty]{}\hat u$ in measure.
\end{proposition}
\begin{proof} Let $\{\hat u_m\}_{m=1}^\infty\subset\Xi(\hat G,\mathcal{Y}(Y))$ be Cauchy. By (\ref{LatticeXiNorm}) and (\ref{Xiu<|u|}) we find that
$$
\|\|\hat u_m(\cdot)\|_\mathcal{Y}-\|\hat u_k(\cdot)\|_\mathcal{Y}\|_\Xi\le c_\Xi\||\|\hat u_m(\cdot)\|_\mathcal{Y}-\|\hat u_k(\cdot)\|_\mathcal{Y}|\|_\Xi
$$
$$
\le c_\Xi K_\Xi\|\|\hat u_m(\cdot)-\hat u_k(\cdot)\|_\mathcal{Y}\|_\Xi=c_\Xi K_\Xi\|\hat u_m-\hat u_k\|_{\mathcal{Y},\Xi},\quad\forall m,k\in\mathbb{N},
$$
which shows that $\{\|\hat u_m(\cdot)\|_\mathcal{Y}\}_{m=1}^\infty$ is Cauchy and hence convergent in $\Xi(\hat G)$,
$$
\|\hat u_m(\cdot)\|_\mathcal{Y}\xrightarrow[m\to\infty]{}\hat w\in\Xi(\hat G).
$$
By the continuous embedding (\ref{XiMeasureProp}), we find that $\{\|\hat u_m(.)\|_\mathcal{Y}\}_{m=1}^\infty$ converges to $\hat w$ in measure as well. On the other hand, by the same embedding (\ref{XiMeasureProp}), the Cauchy sequence $\{\hat u_m\}_{m=1}^\infty\subset\Xi(\hat G,\mathcal{Y}(Y))$ is also Cauchy in $L^0(\hat G,\mathcal{Y}(Y))$. By completeness of $L^0(\hat G,\mathcal{Y}(Y))$ (essentially, a slight modification of Theorem 2.30 in \cite{Fol99}), we find that
$$
\hat u_m\xrightarrow[m\to\infty]{}\hat u\in L^0(\hat G,\mathcal{Y}(Y))
$$
in measure. It follows in particular that $\|\hat u_m(\cdot)\|_\mathcal{Y}$ converges in measure to $\|\hat u(\cdot)\|_\mathcal{Y}$, which shows that $\|\hat u(\cdot)\|_\mathcal{Y}=\hat w\in\Xi(\hat G)$. In view of (\ref{L1locEquiv}) it is clear that $\hat u\in\Xi(\hat G,\mathcal{Y}(Y))$.
\end{proof}

Define the $\mathbb{C}$-conic space
$$
\mathcal{X}(G\times Y)\doteq\mathfrak{F}^{-1}\left(\Xi(\hat G,\mathcal{Y}(Y))\cap\mathcal{D}(\hat G\times Y)'\right)
$$
$$
=\left\{u\in\mathcal{D}(G\times Y)'\,\vline\quad\hat u\in\Xi(\hat G,\mathcal{Y}(Y))\right\}\subset\mathcal{D}(G\times Y)'
$$
and the positive homogeneous function
$$
\|u\|\doteq\|\hat u\|_{\mathcal{Y},\Xi},\quad\forall u\in\mathcal{D}(G\times Y)'.
$$

If the uniform embedding (\ref{UniformEmbed}) holds then by Lemma \ref{EmbeddingLemma}, part 2. we have that
$$
\mathfrak{F}:\mathcal{X}(G\times Y)\to\Xi(\hat G,\mathcal{Y}(Y))
$$
is an isometric bijection and
$$
\mathcal{X}(G\times Y)\hookrightarrow\mathcal{D}(G\times Y)'
$$
is a continuous embedding of topological conic spaces. Again, if $\Xi(\hat G,\mathcal{Y}(Y))\subset L^1_\mathrm{loc}(\hat G,\mathcal{Y}(Y))$ is a vector subspace, then the above maps are linear. And if $\|\cdot\|_{\mathcal{Y},\Xi}$ is a quasi-norm, then the last embedding is uniformly continuous. 
\section{Complex structure and holomorphic functions on $G\times Y$}

Assume now that both $G$ and $Y$ are connected. We need to equip the manifold $G\times Y$ with a complex structure such that $G$ is exactly the ``real'' projection. More precisely, we will assume that
$$
Y\subset\mathbb{R}^n,
$$
so that
$$
G\times Y=G+\imath Y\subset G+\imath\mathbb{R}^n
$$
is an open submanifold with complex coordinates $z=(z_1,\ldots,z_n)\in G\times Y$, where $z_i=x_i+\imath y_i$, $i=1,\ldots,n$. Since $\hat G_i\in\{\mathbb{R},\mathbb{Z}\}$, $i=1,\ldots,n$, the additive bi-character $\langle\cdot,\cdot\rangle$ is well-defined on $Y\times\hat G$, which allows us to extend it analytically to $(G+\imath Y)\times\hat G$,
$$
\langle z,\xi\rangle=\langle x,\xi\rangle+\imath\langle y,\xi\rangle,\quad\chi(z,\xi)=e^{2\pi\imath\langle z,\xi\rangle},\quad\forall z\in G\times Y,\quad\forall\xi\in\hat G.
$$

The following is a variation of the well-known Cauchy-Riemann criterion of holomorphy. For technical convenience here we will work with general, not necessarily $G$-tempered distributions $u\in C^\infty_c(G\times Y)'$, of which the theory is fairly standard. Note that since $C^\infty_c(G\times Y)\hookrightarrow\mathcal{D}(G\times Y)$ continuously, correspondingly, $\mathcal{D}(G\times Y)'\hookrightarrow C^\infty_c(G\times Y)'$ continuously.

\begin{lemma}\label{CRHolLemma} A distribution $u\in C^\infty_c(G\times Y)'$ weakly satisfies the Cauchy-Riemann equations
\begin{equation}
\partial_{\bar z_i}u=0,\quad i=1,\ldots, n,\label{CRWeakEq}
\end{equation}
if and only if $u$ is given by a holomorphic function on $G\times Y$, i.e., $u\in\mathrm{Hol}(G\times Y)$.
\end{lemma}
\begin{proof} The sufficiency part is almost trivial: every holomorphic function $u$ satisfies Cauchy-Riemann equations strongly, \begin{equation}
\partial_{\bar z_i}u(z)=0,\quad\forall z\in G\times Y,\quad i=1,\ldots,n,\label{CRStrongEq}
\end{equation}
and hence also weakly. For the necessity suppose that $u\in C^\infty_c(G\times Y)'$ satisfies (\ref{CRWeakEq}). The principal symbol of the operator $\partial_{\bar z_i}$ is
$$
p_i(z,\zeta)=\frac\imath2(\zeta_{x_i}+\imath\zeta_{y_i}),\quad\forall(z,\zeta)\in T_*(G\times Y)\setminus\{0\},
$$
and the characteristic variety is
$$
\mathrm{Char}(\partial_{\bar z_i})=\left\{(z,\zeta)\in T_*(G\times Y)\setminus\{0\}\,\vline\quad\zeta_{x_i}=\zeta_{y_i}=0\right\}.
$$
By Theorem 8.31 in [H\"ormader I], from (\ref{CRWeakEq}) we have
$$
\mathrm{WF}(u)\subset\bigcap_{i=1}^n\mathrm{Char}(\partial_{\bar z_i})=\emptyset,
$$
whence $u\in C^\infty(G\times Y)$. But then $u$ satisfies the Cauchy-Riemann equations in the strong sense (\ref{CRStrongEq}), and thus, by Hartog's theorem, $u\in\mathrm{Hol}(G\times Y)$.
\end{proof}

It is well-known that $C^\infty_c(\hat G)\subset\mathcal{S}(\hat G)$ is dense, and thus, so is $C^\infty_c(\hat G\times Y)\subset\mathcal{D}(\hat G\times Y)$. It follows that
\begin{align}
\left(\forall\hat u\in C^\infty_c(\hat G\times Y)'\right) & \nonumber\\
\hat u\in\mathcal{D}(\hat G\times Y)'\quad\Leftrightarrow\quad &\left(\forall\{\hat f_m\}_{m=1}^\infty\subset C^\infty_c(\hat G\times Y)\right)\exists\lim_{m\to\infty}\hat f_m\in\mathcal{D}(\hat G\times Y)\,\Rightarrow\,\exists\lim_{m\to\infty}\hat u(f_m).
\end{align}
For every $\hat u\in C^\infty_c(\hat G\times Y)'$ and $\hat\phi\in C^\infty(\hat G\times Y)$, the product $\hat\phi\hat u\in C^\infty_c(\hat G\times Y)'$ is defined by $\hat\phi\hat u(\hat f)=\hat u(\hat\phi\hat f)$ for all $\hat f\in C^\infty_c(\hat G\times Y)$. Similarly, for every $\hat u\in\mathcal{D}(\hat G\times Y)'$ and $\hat\phi\in\mathcal{O}_\mathrm{M}(\hat G\times Y)$, the product $\hat\phi\hat u\in\mathcal{D}(\hat G\times Y)'$ is defined. In both cases the fact is used that $\hat\phi$ acts as a multiplier on the test function space.

Note that the continuous embeddings $\mathcal{S}(\hat G)'\hookrightarrow\mathcal{D}(\hat G\times Y)'$ and $C^\infty_c(\hat G)'\hookrightarrow C^\infty_c(\hat G\times Y)'$ are  given by
$$
\hat u(\hat f)\doteq\int\limits_{Y}\hat u(\hat f(.,y))d\nu(y)
$$
for $\forall\hat f\in\mathcal{D}(\hat G\times Y)$ and $\forall\hat u\in\mathcal{S}(\hat G)'$, or $\forall\hat f\in C^\infty_c(\hat G\times Y)$ and $\forall\hat u\in C^\infty_c(\hat G)'$, respectively.

By $e^{-2\pi\langle\cdot,\cdot\rangle}=\chi(\imath\cdot,\cdot)$ we will denote the function $\hat G\times Y\ni(\xi,y)\mapsto e^{-2\pi\langle y,\xi\rangle}\in\mathbb{R}$. Our next step towards a Fourier characterization of holomorphic function spaces is the following
\begin{proposition}\label{CRFourierProp} The Fourier transform gives a bijection
$$
\bigcap_{i=1}^n\ker\partial_{\bar z_i}\cap\mathcal{D}(G\times Y)'\to\left[e^{-2\pi\langle\cdot,\cdot\rangle}\cdot C^\infty_c(\hat G)'\right]\cap\mathcal{D}(\hat G\times Y)'.
$$
\end{proposition}
\begin{proof} Note that for $i=1,\ldots, n$ we have
$$
\partial_{x_i}\breve f(x,y)=-\breve{[\partial_{x_i}f]}(x,y),\quad\partial_{y_i}\breve f(x,y)=\breve{[\partial_{y_i}f]}(x,y),\quad\forall(x,y)\in G\times Y,\quad\forall f\in\mathcal{D}(G\times Y).
$$
Moreover,
$$
\int\limits_{G\times Y}\breve h(x,y)\partial_{y_i}^\top\breve f(x,y)d\mu(x,y)=\int\limits_{G\times Y}\breve f(x,y)\partial_{y_i}\breve h(x,y)d\mu(x,y)=\int\limits_{G\times Y}\breve f(x,y)\breve{[\partial_{y_i}h]}(x,y)d\mu(x,y)
$$
$$
=\int\limits_{G\times Y}f(x,y)\partial_{y_i}h(x,y)d\mu(x,y)=\int\limits_{G\times Y}h(x,y)\partial_{y_i}^\top f(x,y)d\mu(x,y)=\int\limits_{G\times Y}\breve{[\partial_{y_i}^\top f]}(x,y)\breve h(x,y)d\mu(x,y),
$$
$$
\forall f\in\mathcal{D}(G\times Y),\quad\forall h\in C^\infty_c(G\times Y),\quad i=1,\ldots,n,
$$
where we used the unimodularity of the Abelian group $G$, i.e., $d\mu(-x,y)=d\mu(x,y)$. It follows for $i=1,\ldots,n$ that
$$
\partial_{y_i}^\top\breve f(x,y)=\breve{[\partial_{y_i}^\top f]}(x,y),\quad\forall(x,y)\in G\times Y,\quad\forall f\in\mathcal{D}(G\times Y).
$$
Then for a given $u\in\mathcal{D}(G\times Y)'$ and $i=1,\ldots,n$ the following chain of equivalences holds:
$$
\partial_{\bar z_i}u=0\quad\Leftrightarrow\quad\left(\partial_{x_i}+\imath\partial_{y_i}\right)u=0\quad\Leftrightarrow\quad\left(\forall f\in\mathcal{D}(G\times Y)\right)\quad u\left([-\partial_{x_i}+(\imath\partial_{y_i})^\top]\breve f\right)=0
$$
$$
\Leftrightarrow\quad\left(\forall f\in\mathcal{D}(G\times Y)\right)\quad u\left(\breve{[(\partial_{x_i}+(\imath\partial_{y_i})^\top)f]}\right)=0
$$
$$
\Leftrightarrow\quad\left(\forall f\in\mathcal{D}(G\times Y)\right)\quad\hat u\left(\widehat{[(\partial_{x_i}+(\imath\partial_{y_i})^\top)f]}\right)=0
$$
$$
\Leftrightarrow\quad\left(\forall f\in\mathcal{D}(G\times Y)\right)\quad\hat u\left([2\pi\imath\xi_i+(\imath\partial_{y_i})^\top]\hat f\right)=0
$$
$$
\Leftrightarrow\quad\left(2\pi\imath\xi_i+\imath\partial_{y_i}\right)\hat u=0\quad\Leftrightarrow\quad\partial_{y_i}\left[e^{2\pi\langle\cdot,\cdot\rangle}\cdot\hat u\right]=0.
$$
Here we used (\ref{^fprop}) together with $\partial_{y_i}\in\mathrm{D}(Y)$ and (\ref{D_Y_TransposeEq}). Thus, for $u\in\mathcal{D}(G\times Y)'$,
$$
u\in\bigcap_{i=1}^n\ker\partial_{\bar z_i}\quad\Leftrightarrow\quad\left(\forall i=1,\ldots,n\right)\partial_{\bar z_i}u=0\quad\Leftrightarrow\quad\left(\forall i=1,\ldots,n\right)\partial_{y_i}\left[e^{2\pi\langle\cdot,\cdot\rangle}\cdot\hat u\right]=0
$$
$$
\Leftrightarrow\quad e^{2\pi\langle\cdot,\cdot\rangle}\cdot\hat u\in C^\infty_c(\hat G)'\quad\Leftrightarrow\quad \hat u\in e^{-2\pi\langle\cdot,\cdot\rangle}\cdot C^\infty_c(\hat G)',
$$
which completes the proof.
\end{proof}

\begin{corollary}\label{HolDCorr} Let $u\in\mathrm{Hol}(G\times Y)\cap\mathcal{D}(G\times Y)'$. Then $\exists\hat u_0\in C^\infty_c(\hat G)'$ such that $\hat u=e^{-2\pi\langle\cdot,\cdot\rangle}\,\hat u_0$. The map $u\mapsto\hat u_0$ is a bijection
$$
\mathfrak{F}_0:\mathrm{Hol}(G\times Y)\cap\mathcal{D}(G\times Y)'\to\left\{\hat u_0\in C^\infty_c(\hat G)'\,\vline\quad e^{-2\pi\langle\cdot,\cdot\rangle}\hat u_0\in\mathcal{D}(\hat G\times Y)'\right\}.
$$
\end{corollary}
\begin{proof} Follows immediately from Lemma \ref{CRHolLemma} and Proposition \ref{CRFourierProp}.
\end{proof}

\begin{remark} Remember that $Y$ is connected, and let $0\in Y$. If $u\in\mathrm{Hol}(G\times Y)\cap\mathcal{D}(G\times Y)'$ and $u(\cdot,y)\in L^1(G)$ for $\forall y\in Y$ then we can obtain the above result very easily:
$$
\hat u(\xi,y)e^{2\pi\langle y,\xi\rangle}=\int\limits_Gu(x,y)e^{-2\pi\imath\langle x+\imath y,\xi\rangle}dx=\int\limits_Gu(x,0)e^{-2\pi\imath\langle x,\xi\rangle}dx=\hat u(\xi,0)
$$
$$
=\hat u_0(\xi),\quad\forall(\xi,y)\in\hat G\times Y.
$$
For a general $u\in\mathrm{Hol}(G\times Y)\cap\mathcal{D}(G\times Y)'$, the relation between $u(\cdot,0)$ and $\hat u_0$ is far less obvious.
\end{remark}

\section{Holomorphic functions in the mixed-Fourier-norm spaces}

Denote by $\rho:\hat G\to(0,+\infty]$ the function
$$
\rho(\xi)\doteq\|e^{-2\pi\langle\cdot,\xi\rangle}\|_\mathcal{Y},\quad\forall\xi\in\hat G.
$$
Define the $\mathbb{C}$-conic space
$$
\Xi(\hat G,\rho)\doteq\left\{\hat u\in L^1_\mathrm{loc}(\hat G)\,\vline\quad|\hat u|\rho\in\Xi(\hat G)\right\}.
$$
If $\rho\notin L^0(\hat G)$ then all $\hat u\in\Xi(\hat G,\rho)$ have to vanish on certain subsets in order to make $\hat u\rho\in L^1_\mathrm{loc}(\hat G)$. Introduce the positive homogeneous function
$$
\|\hat u\|_{\Xi,\rho}\doteq\||\hat u|\rho\|_\Xi,\quad\forall\hat u\in L^1_\mathrm{loc}(\hat G)\,\cap\,\frac1{\rho}\cdot L^1_\mathrm{loc}(\hat G).
$$
We equip $\Xi(\hat G,\rho)$ with the coarsest topology which makes $\|\cdot-\hat v\|_{\Xi,\rho}$ continuous for $\forall v\in\Xi(\hat G,\rho)$. This demonstrates $\Xi(\hat G,\rho)$ as a topological conic space. Setting $\|\hat u\|_{\Xi,\rho}\doteq+\infty$ for all $\hat u\in C^\infty_c(\hat G)'$ such that $\hat u\notin L^1_\mathrm{loc}(\hat G)$ or $\hat u\cdot\rho\not\in L^1_\mathrm{loc}(\hat G)$ we can assume without loss of generality that $\|\hat u\|_{\Xi,\rho}$ is defined for all $\hat u\in C^\infty_c(\hat G)'$.

\begin{lemma}\label{LatticeLemmaScalar} The following statements hold:
\begin{itemize}

\item[1.] If the lattice property (\ref{LatticeXiSpace}) is true then $\Xi(\hat G,\rho)\subset L^1_\mathrm{loc}(\hat G)$ is a vector subspace.

\item[2.] If the quasi-lattice property (\ref{LatticeXiNorm}) is true then $\|\cdot\|_{\Xi,\rho}$ is a quasi-norm with the same triangle constant $K_\Xi$.

\end{itemize}
\end{lemma}
\begin{proof} The proof is virtually identical to that of Lemma \ref{LatticeLemma}.
\end{proof}

If $\|\cdot\|_{\Xi,\rho}$ is a quasi-norm then $\Xi(\hat G,\rho)$ acquires the structure of a uniform conic space, with a base of entourages
$$
\left\{(\hat u,\hat v)\in\Xi(\hat G,\rho)^2\,\vline\quad\|\hat u-\hat v\|_{\Xi,\rho}<\epsilon\right\}.
$$

\begin{proposition}\label{XirhoCompleteProp} Suppose that $\|\cdot\|_{\Xi,\rho}$ is a quasi-norm, and conditions (\ref{Xi|u|<u}), (\ref{Xiu<|u|}), (\ref{u_|u|}) and (\ref{|u|_u}) hold true. If $\frac1\rho\in L^\infty_\mathrm{loc}(\hat G)$ then the uniform space $\Xi(\hat G,\rho)$ is complete.
\end{proposition}
\begin{proof} Let $\{\hat u_m\}_{m=1}^\infty\subset\Xi(\hat G,\rho)$ be Cauchy, i.e.,
$$
\left(\forall\epsilon>0\right)\left(\exists N_\epsilon\in\mathbb{N}\right)\left(\forall m,k>N_\epsilon\right)\quad\|\hat u_m-\hat u_k\|_{\Xi,\rho}<\epsilon.
$$
From $|\hat u_m|\rho\in\Xi(\hat G)$ and (\ref{|u|_u}) we have $\hat u_m\cdot\rho\in\Xi(\hat G)$ for $\forall m\in\mathbb{N}$. Then, in view of (\ref{Xiu<|u|}), from
$$
\|\hat u_m\cdot\rho-\hat u_k\cdot\rho\|_\Xi\le c_\Xi\||\hat u_m-\hat u_k|\rho\|_\Xi=c_\Xi\|\hat u_m-\hat u_k\|_{\Xi,\rho}
$$
we find that $\{\hat u_m\cdot\rho\}_{m=1}^\infty\subset\Xi(\hat G)$ is Cauchy and hence $\hat u_m\cdot\rho\to\hat w\in\Xi(\hat G)$. From $\hat w\in L^1_\mathrm{loc}(\hat G)$ and $\frac1\rho\in L^\infty_\mathrm{loc}(\hat G)$ it follows that $\hat u\doteq\frac{\hat w}\rho\in L^1_\mathrm{loc}(\hat G)$. Moreover, $\hat w\in\Xi(\hat G)$ implies by (\ref{u_|u|}) that $|\hat u|\rho=|\hat w|\in\Xi(\hat G)$, so that $\hat u\in\Xi(\hat G,\rho)$. Finally, using (\ref{Xi|u|<u}) we see that
$$
\|\hat u_m-\hat u\|_{\Xi,\rho}=\||\hat u_m-\hat u|\rho\|_\Xi=\||\hat u_m\cdot\rho-\hat w|\|_\Xi\le C_\Xi\|\hat u_m\cdot\rho-\hat w\|_\Xi\xrightarrow[m\to\infty]{}0,
$$
which shows that $\hat u_m\to\hat u$ in $\Xi(\hat G,\rho)$.
\end{proof}

The following unsurprising lemma is included for completeness.

\begin{lemma}\label{hatu_0L1locLemma} Suppose that the uniform embedding (\ref{UniformEmbed}) is true, so that, by Lemma \ref{EmbeddingLemma}, the embedding (\ref{L1locEmbed}) holds. Then for every $\hat u_0\in C^\infty_c(\hat G)'$,
$$
e^{-2\pi\langle\cdot,\cdot\rangle}\hat u_0\in L^1_\mathrm{loc}(\hat G,\mathcal{Y}(Y))\quad\Rightarrow\quad\hat u_0\in L^1_\mathrm{loc}(\hat G).
$$
\end{lemma}
\begin{proof} Fix $h\in C^\infty_c(Y)$ with $h\ge0$, and define $\hat F_h\in C^\infty(\hat G,\mathbb{R}_+)$ by
$$
F_h(\xi)\doteq\int\limits_Ye^{-2\pi\langle y,\xi\rangle}h(y)d\nu(y),\quad\forall\xi\in\hat G.
$$
Denote $\hat w\doteq e^{-2\pi\langle\cdot,\cdot\rangle}\hat u_0\in L^1_\mathrm{loc}(\hat G,\mathcal{Y}(Y))$ and define $\hat H_h:\hat G\to\mathbb{C}$ by
$$
\hat H_h(\xi)\doteq\hat w(\xi)(h),\quad\forall\xi\in\hat G.
$$
That $\hat H_h$ is Borel-measurable follows from the strong, and hence also weak measurability of $\hat w:\hat G\to\mathcal{Y}(Y)$. By (\ref{UniformEmbed}),
$$
|\hat H_h(\xi)|\le\wp(h)\|\hat w(\xi)\|_\mathcal{Y},\quad\forall\xi\in\hat G,
$$
and since $\|\hat w(\cdot)\|_\mathcal{Y}\in L^1_\mathrm{loc}(\hat G)$ by (\ref{L1locEquiv}), we have that $\hat H_h\in L^1_\mathrm{loc}(\hat G)$. Since $(\hat F_h)^{-1}=1/\hat F_h\in C^\infty(\hat G,\mathbb{R}_+)$, we can write
$$
\hat u_0(\hat f)=\hat u_0\left((\hat F_h)^{-1}\hat f\cdot\int\limits_Ye^{-2\pi\langle y,\xi\rangle}h(y)d\nu(y)\right)=\left[e^{-2\pi\langle\cdot,\cdot\rangle}\hat u_0\right]\left((\hat F_h)^{-1}\hat fh\right)
$$
$$
=\int\limits_{\hat G}\hat F_h(\xi)^{-1}\hat H_h(\xi)\hat f(\xi)d\xi,\quad\forall\hat f\in C^\infty_c(\hat G),
$$
where we used the embedding (\ref{L1locEmbed}). It follows that
$$
\hat u_0=(\hat F_h)^{-1}\cdot\hat H_h\in L^1_\mathrm{loc}(\hat G),
$$
as desired.
\end{proof}

Introduce the $\mathbb{C}$-conic subspace
$$
\mathcal{A}_\mathcal{X}(G\times Y)\doteq\mathcal{X}(G\times Y)\cap\mathrm{Hol}(G\times Y)\subset\mathcal{X}(G\times Y)
$$
consisting of only holomorphic elements.

\begin{proposition}\label{MainProp} Suppose that the uniform embedding (\ref{UniformEmbed}) holds.
\begin{itemize}

\item[1.] The map $u\mapsto\hat u_0$ from Corollary \ref{HolDCorr} restricts to an isometry
\begin{equation}
\mathfrak{F}_0:\mathcal{A}_\mathcal{X}(G\times Y)\to\Xi(\hat G,\rho).\label{F_0A_X}
\end{equation}

\item[2.] If the map
$$
\hat G\ni\xi\mapsto \frac{e^{-2\pi\langle\cdot,\xi\rangle}}{\rho(\xi)}\in C^\infty_c(Y)'
$$
equals almost everywhere a Bochner-measurable function $\hat G\to\mathcal{Y}(Y)$, then the isometry (\ref{F_0A_X}) is surjective.

\end{itemize}
\end{proposition}
\begin{proof} We begin by proving the statement 1. Let $u\in\mathcal{A}_\mathcal{X}(G\times Y)\subset\mathcal{D}(G\times Y)'\cap\mathrm{Hol}(G\times Y)$. Then by Corollary \ref{HolDCorr} we have
$$
\hat u_0\in C^\infty_c(\hat G)',\quad\hat u=e^{-2\pi\langle\cdot,\cdot\rangle}\hat u_0\in\mathcal{D}(\hat G\times Y)'.
$$
But also $u\in\mathcal{X}(G\times Y)$, therefore
$$
\hat u=e^{-2\pi\langle\cdot,\cdot\rangle}\hat u_0\in\Xi(\hat G,\mathcal{Y}(Y)),
$$
whence by (\ref{L1locEquiv}),
$$
e^{-2\pi\langle\cdot,\cdot\rangle}\hat u_0\in L^1_\mathrm{loc}(\hat G,\mathcal{Y}(Y)),\quad\|[e^{-2\pi\langle\cdot,\cdot\rangle}\hat u_0](\cdot)\|_\mathcal{Y}\in\Xi(\hat G).
$$
By Lemma \ref{hatu_0L1locLemma}, $\hat u_0\in L^1_\mathrm{loc}(\hat G)$ and
$$
\|[e^{-2\pi\langle\cdot,\cdot\rangle}\hat u_0](\cdot)\|_\mathcal{Y}=|\hat u_0|\rho\in\Xi(\hat G).
$$
Moreover,
$$
\|u\|=\|\hat u\|_{\mathcal{Y},\Xi}=\left\|e^{-2\pi\langle\cdot,\cdot\rangle}\hat u_0\right\|_{\mathcal{Y},\Xi}=\left\|\left\|[e^{-2\pi\langle\cdot,\cdot\rangle}\hat u_0](\cdot)\right\|_\mathcal{Y}\right\|_\Xi=\||\hat u_0|\rho\|_\Xi=\|\hat u_0\|_{\Xi,\rho},
$$
showing that the map (\ref{F_0A_X}) is an isometry.

For the statement 2., take $\hat u_0\in\Xi(\hat G,\rho)$, i.e., $\hat u_0\in L^1_\mathrm{loc}(\hat G)$ with $|\hat u_0|\rho\in\Xi(\hat G)$. Then $\hat u_0\rho\in L^1_\mathrm{loc}(\hat G)$, and
$$
e^{-2\pi\langle\cdot,\cdot\rangle}\hat u_0=\frac{e^{-2\pi\langle\cdot,\cdot\rangle}}{\rho}\hat u_0\rho
$$
equals almost everywhere a Bochner-measurable function. From
$$
\|[e^{-2\pi\langle\cdot,\cdot\rangle}\hat u_0](\cdot)\|_\mathcal{Y}=|\hat u_0|\rho\in\Xi(\hat G)\subset L^1_\mathrm{loc}(\hat G)
$$
we see that, by Lemma \ref{EmbeddingLemma},
$$
e^{-2\pi\langle\cdot,\cdot\rangle}\hat u_0\in\Xi(\hat G,\mathcal{Y}(Y))\subset\mathcal{D}(\hat G\times Y)',
$$
and thus,
$$
u\doteq\mathfrak{F}^{-1}\left[e^{-2\pi\langle\cdot,\cdot\rangle}\hat u_0\right]\in\mathcal{X}(G\times Y).
$$
But by Corollary \ref{HolDCorr} we also have that $u\in\mathcal{D}(G\times Y)'\cap\mathrm{Hol}(G\times Y)$, so that finally $u\in\mathcal{A}_\mathcal{X}(G\times Y)$ and $\hat u_0=\mathfrak{F}_0u$. This completes the proof that the map (\ref{F_0A_X}) is surjective, and hence an isometric bijection.
\end{proof}

\begin{remark} If
$$
\Xi(\hat G,\mathcal{Y}(Y))\subset L^1_\mathrm{loc}(\hat G,\mathcal{Y}(Y))
$$
is a vector subspace then so is
$$
\mathcal{A}_\mathcal{X}(G\times Y)\subset\mathrm{Hol}(G\times Y)\cap\mathcal{D}(G\times Y)'.
$$
In that case the image
$$
\mathfrak{F}_0\left(\mathcal{A}_\mathcal{X}(G\times Y)\right)\subset C^\infty_c(\hat G)'
$$
is a vector subspace, which coincides with $\Xi(\hat G,\rho)$ if the map (\ref{F_0A_X}) is surjective.
\end{remark}

\begin{remark} Suppose that the quasi-lattice property (\ref{LatticeXiNorm}) is true, so that by, Lemma \ref{LatticeLemma} and Lemma \ref{LatticeLemmaScalar}, both $\mathcal{A}_\mathcal{X}(G\times Y)$ and $\Xi(\hat G,\rho)$ are uniform spaces. If $\Xi(\hat G,\rho)$ is complete and if the isometric map (\ref{F_0A_X}) is surjective, then $\mathcal{A}_\mathcal{X}(G\times Y)$ is also complete.
\end{remark}

\section{Analytic factorisation}

The entire analysis of the preceding sections was based on the factorised space $G\times Y$ - the trivial Abelian principal bundle. However, already in the basic, archetypical examples one does not begin with such a domain initially. Instead, one begins with an analytical domain $\Omega$ with an action of an Abelian Lie group $G$. Before the analysis carried out above can be applied, one has to recognise the principal bundle $G\times Y$ inside $\Omega$, or to analytically factorise $\Omega$. The present section is devoted to the subtleties arising in that process.

Let us refer back to Section \ref{MixedFourierNormSection} for notations and definitions. Consider the monoid
$$
\overline{\hat G_+}\doteq\prod_{i=1}^n\left(\hat G_i\cap[0,+\infty)\right).
$$
Since $G$ is connected, it is clear that $\hat G_i\cap[0,+\infty)\in\{[0,+\infty),\mathbb{N}_0\}$ for $i=1,\ldots,n$. In what follows we will always assume that
\begin{equation}
Y\subset(\mathbb{R}_+)^n.\label{y>0}
\end{equation}
This is the analogue of $y>0$ for the upper half-plane or $r<1$ for the unit disk, as models for a good domain of analytical functions. Below we will see that, with respect to the chosen complex structure, $y_i\to 0$ approaches an unrestricted edge of the domain, while at $y_i\to+\infty$ we expect bounded behaviour, $i=1,\ldots,n$.

The main motivation behind the introduction of $\overline{\hat G_+}$ and the restriction (\ref{y>0}) is the following Paley-Wiener-type property expected to hold in most common situations:
\begin{equation}
\supp\hat u_0\in\overline{\hat G_+},\quad\forall u\in\mathcal{A}_\mathcal{X}(G\times Y).\label{PW}
\end{equation}
Below is a sufficient condition for the property (\ref{PW}).

\begin{proposition}\label{PWSuppProp} Suppose that the uniform embedding (\ref{UniformEmbed}) is true. If
$$
\rho\;\vline_{\,\overline{\hat G_+}^\complement}\equiv+\infty,
$$
then the property (\ref{PW}) holds.
\end{proposition}
\begin{proof} Let (\ref{UniformEmbed}) be true and take $\forall u\in\mathcal{A}_\mathcal{X}(G\times Y)$. By Proposition \ref{MainProp}, we have $\hat u_0\in\Xi(\hat G,\rho)$. This implies $|\hat u_0|\rho\in\Xi(\hat G)\subset L^1_\mathrm{loc}(\hat G)$, whence
$$
\supp\hat u_0\subset\supp\frac1\rho\subset\overline{\hat G_+}.
$$
The proof is complete.
\end{proof}

For $y_1,y_2\in\mathbb{R}^n$, we will write $y_1>y_2$ if $y_1-y_2\in(\mathbb{R}_+)^n$. The property (\ref{PW}) is the precursor of the desired boundedness at $y_i\to+\infty$, $i=1,\ldots,n$:
\begin{equation}
(\exists y_0>0)\quad(\forall u\in\mathcal{A}_\mathcal{X}(G\times Y))\quad\sup_{y>y_0}\|u(\cdot,y)\|_\infty<\infty.\label{BDProp}
\end{equation}

\begin{proposition}\label{PWBddProp} Suppose that the uniform embedding (\ref{UniformEmbed}) and the property (\ref{PW}) hold, and $\exists y_0\ge0$ such that
$$
\frac{e^{-2\pi\langle y_0,\cdot\rangle}}{\rho}\vline_{\,\overline{\hat G_+}}\in L^\infty(\,\overline{\hat G_+}\,).
$$
Then the boundedness property (\ref{BDProp}) is valid.
\end{proposition}
\begin{proof} Let $u\in\mathcal{A}_\mathcal{X}(G\times Y)$, so that, by Corollary \ref{HolDCorr}, $\hat u=\hat u_0e^{-2\pi\langle\cdot,\cdot\rangle}$, where, by the statement 1. of Proposition \ref{MainProp}, $\hat u_0\in\Xi(\hat G,\rho)$. In particular, $\hat u_0\in L^1_\mathrm{loc}(\hat G)$ and $|\hat u_0|\rho\in\Xi(\hat G)\subset\mathcal{S}(\hat G)'$. Then for $\forall y>y_1>y_0$ we have
$$
\|u(\cdot,y)\|_\infty\le\int\limits_{\overline{\hat G_+}}|\hat u(\xi,y)|d\xi=\int\limits_{\overline{\hat G_+}}|\hat u_0(\xi)|\rho(\xi)\frac{e^{-2\pi\langle y,\xi\rangle}}{\rho(\xi)}d\xi
$$
$$
\le\int\limits_{\overline{\hat G_+}}|\hat u_0(\xi)|\rho(\xi)\frac{e^{-2\pi\langle y_0,\xi\rangle}}{\rho(\xi)}e^{-2\pi\langle y_1-y_0,\xi\rangle}d\xi\le\left\|\frac{e^{-2\pi\langle y_0,\cdot\rangle}}{\rho}\right\|_{L^\infty(\overline{\hat G_+})}\cdot[|\hat u_0|\rho](f_\epsilon)<\infty,
$$
where $f_\epsilon\in\mathcal{S}(\hat G)$ with $f|_{\overline{\hat G_+}}=e^{-2\pi\langle y_1-y_0,\cdot\rangle}$. The proof is complete.
\end{proof}

Finally, let $\Omega$ be a complex manifold and $\Phi:G\times Y\to\Omega$ a bi-holomorphism onto an open subset $\Phi(G\times Y)\subset\Omega$. It is clear that a function $u:\Omega\to\mathbb{C}$ is holomorphic in $\Phi(G\times Y)$ if and only if $\Phi^*u\in\mathrm{Hol}(G\times Y)$. Moreover,

\begin{remark} If the boundedness condition (\ref{BDProp}) is true, then
$$
\Phi^*u\in\mathcal{A}_\mathcal{X}(G\times Y)\quad\Rightarrow\quad\sup\left\{|u(z)|\;\vline\quad z\in\overline{\Phi\bigl(G\times(Y\cap[y_0+(\mathbb{R}_+)^n])\bigr)}\right\}<\infty.
$$
\end{remark}

\section{Examples}

\subsection{Mixed-Fourier-norm space: the elliptic geometry case}\label{section-1}

Our first example, which has served as the main motivation for this study, is the unit disc with the elliptic model. Namely, we begin with the domain $\Omega=\mathbb{D}\subset\mathbb{C}$ and the Lie group $\mathbb{T}$ acting as $w\mapsto e^{2\pi\imath x}w$ on $w\in\Omega$. On the open dense subset $\mathbb{D}\setminus\{0\}$, this action is free. Our model space is $G\times Y=\mathbb{T}\times\mathbb{R}_+$, which is mapped to $\mathbb{D}\setminus\{0\}$ by the $G$-equivariant biholomorphism $\Phi:\mathbb{T}\times\mathbb{R}_+\to\mathbb{D}$ given as
$$
\Phi(z)=\Phi(x+\imath y)=e^{2\pi z}=e^{2\pi(\imath x-y)},\quad\forall z\in\mathbb{T}\times\mathbb{R}_+.
$$
This mapping has the important property that $w\to\partial\mathbb{D}$ corresponds to $z\to\mathbb{T}\times\{0\}$, while $w\to0$ corresponds to $z\to+\imath\infty$. The normalized Lebesgue measure on $\mathbb{D}$ corresponds to the measure $dxd\nu(y)$ on $\mathbb{T}\times\mathbb{R}_+$, where
$$
d\nu(y)=4\pi e^{-4\pi y}dy.
$$
The compactness of the group $G$ provides remarkable technical advantages compared to the general case. In particular, all distributions are $G$-tempered,
$$
\mathcal{D}(G\times Y)=C^\infty_c(\mathbb{T}\times\mathbb{R}_+),\quad\mathcal{D}(G\times Y)'=C^\infty_c(\mathbb{T}\times\mathbb{R}_+)'.
$$
Another advantage of a compact $G$ is the discreteness of $\hat G=\mathbb{Z}$. Let us take $\Xi(\hat G)\doteq l^q$ for some $q\in[1,+\infty)$, and $\mathcal{Y}(Y)\doteq X(\mathbb{R}_+)$ a Banach space of distributions uniformly embedded in $C^\infty_c(\mathbb{R}_+)'$ in the sense of (\ref{UniformEmbed}). The first important object is the space
$$
\Xi(\hat G,\mathcal{Y}(Y))=\left\{\hat u:\mathbb{Z}\to X(\mathbb{R}_+)\;\vline\quad\|\hat u(\cdot)\|_{X(\mathbb{R}_+)}\in l^q\right\},
$$
which can also be written as
$$
\Xi(\hat G,\mathcal{Y}(Y))=\left\{\{\hat u_\xi\}_{\xi=-\infty}^\infty\subset X(\mathbb{R}_+)\;\vline\quad\{\|\hat u_\xi\|_{X(\mathbb{R}_+)}\}_{\xi=-\infty}^\infty\subset l^q\right\}.
$$
It is easy to check that the (quasi-)lattice properties (\ref{LatticeXiNorm}) and (\ref{LatticeXiSpace}) are satisfied, therefore, by Lemma \ref{LatticeLemma}, we have that $\Xi(\hat G,\mathcal{Y}(Y))$ is a normed vector space with the norm
$$
\|\hat u\|_{\mathcal{Y},\Xi}=\|\|\hat u(\cdot)\|_{X(\mathbb{R}_+)}\|_q=\left(\sum_{\xi=-\infty}^\infty\|\hat u_\xi\|_{X(\mathbb{R}_+)}^q\right)^{\frac1q}.
$$
By Lemma \ref{EmbeddingLemma}, we have a continuous linear embedding
$$
\Xi(\hat G,\mathcal{Y}(Y))\hookrightarrow\mathcal{D}(\hat G\times Y)'=C^\infty_c(\mathbb{Z}\times\mathbb{R}_+)'.
$$
Since the condition (\ref{Xiu<|u|}) is trivially true, and $l^q\hookrightarrow l^0$ continuously, we have by Proposition \ref{XiMeasureClosed} that Cauchy sequences in $\Xi(\hat G,\mathcal{Y}(Y))$ converge at least in measure.

The mixed-Fourier-norm space $\mathcal{X}(G\times Y)$ is introduced as
$$
\mathcal{X}(G\times Y)=\left\{u\in C^\infty_c(\mathbb{T}\times\mathbb{R}_+)'\;\vline\quad\hat u\in\Xi(\hat G,\mathcal{Y}(Y))\right\},
$$
with the norm $\|u\|=\|\hat u\|_{\mathcal{Y},\Xi}$. It follows immediately that
$$
\mathfrak{F}:\mathcal{X}(G\times Y)\to\Xi(\hat G,\mathcal{Y}(Y))
$$
is an isometric isomorphism of normed vector spaces.

The function (sequence) $\rho:\mathbb{Z}\to(0,+\infty]$ is defined as
$$
\rho(\xi)=\left\|e^{-2\pi\langle\cdot,\xi\rangle}\right\|_{X(\mathbb{R}_+)},\quad\forall\xi\in\mathbb{Z}.
$$
It is clear that $\frac1\rho\in L^\infty_\mathrm{loc}(\hat G)=l^0$. The weighted sequence space $\Xi(\hat G,\rho)$ is given as
$$
\Xi(\hat G,\rho)=\left\{\hat u:\mathbb{Z}\to\mathbb{C}\;\vline\quad \hat u\rho\in l^q\right\}.
$$
with the norm $\|\hat u\|_{\Xi,\rho}=\|\hat u\rho\|_q$. By Lemma \ref{LatticeLemmaScalar}, $\Xi(\hat G,\rho)$ is a normed vector space (which was already clear at this point). What is more important, since properties (\ref{Xi|u|<u}), (\ref{Xiu<|u|}), (\ref{u_|u|}) and (\ref{|u|_u}) are obviously true, by Proposition \ref{XirhoCompleteProp}, the normed vector space $\Xi(\hat G,\rho)$ is actually complete, i.e., a Banach space.

The holomorphic mixed-Fourier-norm space $\mathcal{A}_\mathcal{X}(G\times Y)$ consists of all members of the space $\mathcal{X}(G\times Y)$ which are holomorphic in $\mathbb{T}\times\mathbb{R}_+$,
$$
\mathcal{A}_\mathcal{X}(G\times Y)=\mathcal{X}(G\times Y)\,\cap\,\mathrm{Hol}(G\times Y)
$$
$$
=\left\{u\in\mathrm{Hol}(\mathbb{T}\times\mathbb{R}_+)\;\vline\quad\hat u\in\Xi(\hat G,\mathcal{Y}(Y))\right\}.
$$
Since the discrete function (sequence)
$$
\mathbb{Z}\ni\xi\mapsto\frac{e^{-2\pi\langle\cdot,\xi\rangle}}{\rho(\xi)}\in X(\mathbb{R}_+)
$$
is obviously measurable (with the convention $\frac1{+\infty}=0$), by Proposition \ref{MainProp} we have that
$$
\mathfrak{F}:\mathcal{A}_\mathcal{X}(G\times Y)\to e^{-2\pi\langle\cdot,\cdot\rangle}\,\cdot\,\Xi(\hat G,\rho),\quad u(z)\mapsto e^{-2\pi\xi y}\,\hat u_0(\xi),
$$
is an isometric isomorphism of normed vector spaces. In particular, since $\Xi(\hat G,\rho)$ is complete, so is $\mathcal{A}_\mathcal{X}(G\times Y)$.

Now suppose that the conditions of Proposition \ref{PWSuppProp} are satisfied. This gives us a Paley-Wiener-type support property,
\begin{equation}
\mathrm{supp}\,\hat u_0\subset\mathbb{N}_0,\quad\forall u\in\mathcal{A}_\mathcal{X}(G\times Y).\label{PWSuppE}
\end{equation}
Moreover, if the sequence
\begin{equation}
\frac{e^{-2\pi\langle y_0,\cdot\rangle}}{\rho}\in l^\infty\label{Cond1}
\end{equation}
for some $y_0\ge0$, then by Proposition \ref{PWBddProp} we have a Paley-Wiener-type boundedness property,
\begin{equation}
\sup_{y>y_1}\|u(\cdot,y)\|_\infty<\infty,\quad\forall u\in\mathcal{A}_\mathcal{X}(G\times Y),\label{PWBddE}
\end{equation}
for every $y_1>y_0$.

Now let us translate all this back to the original domain $\Omega=\mathbb{D}$. We have
$$
\Phi^*C^\infty_c(\mathbb{D}\setminus\{0\})=C^\infty_c(\mathbb{T}\times\mathbb{R}_+),\quad\Phi^*C^\infty_c(\mathbb{D}\setminus\{0\})'=C^\infty_c(\mathbb{T}\times\mathbb{R}_+)'.
$$
For $f\in C^\infty_c(\mathbb{D}\setminus\{0\})$, its Fourier transform is the sequence of distributions $\{\hat f_\xi\}_{\xi=-\infty}^\infty\subset C^\infty_c((0,1))'$ given by (formally)
$$
\hat f_\xi(r)=(\mathfrak{F}\Phi^*f)_\xi(r)=\int\limits_\mathbb{T}f(re^{2\pi\imath x})e^{-2\pi\imath x\xi}dx,
$$
where $r(y)=e^{-2\pi y}$. This map $r:\mathbb{R}_+\to(0,1)$ maps an arbitrary uniformly embedded Banach space $X((0,1))\hookrightarrow C^\infty_c((0,1))'$ to an arbitrary uniformly embedded Banach space
$$
r^*X((0,1))=X(\mathbb{R}_+)\hookrightarrow C^\infty_c(\mathbb{R}_+).
$$
If $u\in C^\infty_c(\mathbb{T}\times\mathbb{R}_+)'$ and $f\in C^\infty_c(\mathbb{D}\setminus\{0\})'$ such that $u=\Phi^*f$, then
$$
\hat u_\xi=r^*\hat f_\xi,\quad\hat u_\xi\in X(\mathbb{R}_+)\quad\Leftrightarrow\quad\hat f_\xi\in X((0,1)),\quad\forall\xi\in\mathbb{Z}.
$$
The mixed-Fourier-norm space $\mathcal{L}^{q;X}(\mathbb{D})$ is defined to be the pushforward of $\mathcal{X}(G\times Y)$,
$$
\Phi^*\mathcal{L}^{q;X}(\mathbb{D})=\mathcal{X}(G\times Y).
$$
More precisely,
$$
\mathcal{L}^{q;X}(\mathbb{D})=\left\{f\in C^\infty_c(\mathbb{D}\setminus\{0\})'\;\vline\quad\{\hat f_\xi\}_{\xi=-\infty}^\infty\subset X((0,1)),\quad\|f\|<\infty\right\},
$$
with the norm
$$
\|f\|=\left(\sum_{\xi=-\infty}^\infty\|\hat f_\xi\|_{X((0,1))}^q\right)^{\frac1q}.
$$
If we assume in addition that $X((0,1))\subset L^1((0,1))$ then $f\in\mathcal{L}^{q;X}(\mathbb{D})$ can be extended uniquely to $f\in C^\infty_c(\mathbb{D})'$. However, note that
$$
\mathcal{L}^{q;X}(\mathbb{D})\neq\left\{f\in C^\infty_c(\mathbb{D})\;\vline\quad\{\hat f_\xi\}_{\xi=-\infty}^\infty\subset X((0,1)),\quad\|f\|<\infty\right\}.
$$
Indeed, take $f=\delta$ the Dirac distribution. Then $\hat f=0$, which destroys the isomorphism with $\Xi(\hat G,\mathcal{Y}(Y))$. Only those distributions on the disc which have no singular support at $0$ can be taken in $\mathcal{L}^{q;X}(\mathbb{D})$, and this cannot be checked on the Fourier side.

Recall that we were only able to show on general grounds that $\Xi(\hat G,\mathcal{Y}(Y))$ is complete in the topology of convergence in measure. However, it was shown (Theorem 3.1 in \cite{KS-2017}) that, if all step functions are in $X((0,1))$, then the space $\mathcal{L}^{q;X}(\mathbb{D})$ is complete.

Finally, the mixed-Fourier-norm Bergman space $\mathcal{A}^{q;X}(\mathbb{D})$ is defined as the pushforward of the holomorphic mixed-Fourier-norm space $\mathcal{A}_\mathcal{X}(G\times Y)$,
$$
\Phi^*\mathcal{A}^{q;X}(\mathbb{D})=\mathcal{A}_\mathcal{X}(G\times Y).
$$
In precise terms,
$$
\mathcal{A}^{q;X}(\mathbb{D})=\left\{f\in\mathrm{Hol}(\mathbb{D}\setminus\{0\})\;\vline\quad\{\hat f_\xi\}_{\xi=-\infty}^\infty\subset X((0,1)),\quad\|f\|<\infty\right\}.
$$
If the Banach space $X((0,1))$ has the property that
$$
\|r^{-\xi}\|_{X((0,1))}=+\infty,\quad\forall\xi\in\mathbb{N},
$$
which is true if $X((0,1))\subset L^1((0,1))$, then we have the support property (\ref{PWSuppE}), giving
$$
\|f\|=\left(\sum_{\xi=0}^\infty\|\hat f_\xi\|_{X((0,1))}^q\right)^{\frac1q}.
$$
On the other hand, if the Banach space $X(\mathbb{R}_+)=\Phi^*X((0,1))$ satisfies the condition (\ref{Cond1}) above, which translates to
$$
\left\{\frac{r_0^\xi}{\|r^\xi\|_{X((0,1))}}\right\}_{\xi=0}^\infty\subset l^\infty
$$
for some $r_0\in(0,1)$, then we have the boundedness result (\ref{PWBddE}), which translates to
$$
\sup_{r<r_0,\,x\in\mathbb{T}}|f(re^{2\pi\imath x})|<\infty,\quad\forall f\in\mathcal{A}^{q;X}(\mathbb{D}).
$$
But then by Hartog's extension theorem, every $f\in\mathcal{A}^{q;X}(\mathbb{D})$ can be uniquely extended to $f\in\mathrm{Hol}(\mathbb{D})$. Thus,
$$
\mathcal{A}^{q;X}(\mathbb{D})=\left\{f\in\mathrm{Hol}(\mathbb{D})\;\vline\quad\{\hat f_\xi\}_{\xi=0}^\infty\subset X((0,1)),\quad\|f\|<\infty\right\}.
$$
Recall that the space $\mathcal{A}_\mathcal{X}(G\times Y)$ was shown to be complete, which implies the completeness of the mixed-Fourier-norm Bergman space $\mathcal{A}^{q;X}(\mathbb{D})$. Note that the completeness of $\mathcal{A}^{q;X}(\mathbb{D})$ was obtained in \cite{KS-2017} under additional assumptions that
all step functions are in $X((0,1))$, and
\begin{equation}\label{eq:COND} n \left|\int_I \tau^n g(\tau) 2\tau d\tau \right|\|r^n\|_{X((0,1))}\leqslant C_0 \|g\|_{X((0,1))},\ n\to\infty, \ \forall g\in
X((0,1)).
\end{equation}
Most importantly, Fourier transform $f\mapsto\hat f$ is an isometric isomorphism of Banach spaces
$$
\mathcal{A}^{q;X}(\mathbb{D})\to\left\{\{\hat f_\xi\, r^\xi\}_{\xi=0}^\infty\;\vline\quad\hat f\in l^q_+(\rho^q)\right\},
$$
$$
l^q_+(\rho^q)=\left\{\hat f\in l^q_+\;\vline\quad\|\hat f\|_{q,\rho}<\infty\right\},\quad\|\hat f\|_{q,\rho}=\left(\sum_{\xi=0}^\infty|f_\xi|^q\|r^\xi\|_{X((0,1))}^q\right)^{\frac1q}.
$$

\begin{remark}
As already noted in the Introduction, the general definition of spaces with a mixed-Fourier-norm on the unit disc associated with polar coordinates was given in the paper \cite{KS-2017} within the framework of classical harmonic analysis. Further, special cases of $X$ spaces were investigated in the works \cite{KS-2018,KS-2016-1,KS-2016-2,KSmi-CVEE,KSmi-IzvVuz,K-Grand-2021}. In particular, results of Paley-Wiener and other types are presented in the papers \cite{KSmi-CVEE,KSmi-IzvVuz}. The idea of this Section is first of all to apply our general group-based approach to this situation and to clarify our statements in this important particular case. This gives us a holistic and deeper look at the problem.
\end{remark}

\subsection{Mixed-Fourier-norm space: the parabolic geometry case}\label{section-3}

We will realise the parabolic geometry on the upper half-plane $\Pi=\mathbb{R}\times\mathbb{R}_+$. In this case, we are already on the model space $G\times Y=\mathbb{R}\times\mathbb{R}_+$, with $\mathbb{R}$ acting naturally as $z\mapsto x+z$ on $z\in\Pi$, and there is no need for a holomorphic transformation. In this case our measure on $\mathbb{R}_+$ will be
$$
d\nu_\lambda(y)=(\lambda+1)(2y)^\lambda dy
$$
for some $\lambda>-1$, which will be fixed for the rest of this section. Our test function and distribution spaces will be
$$
\mathcal{D}(G\times Y)=\mathcal{S}(\mathbb{R})\,\hat{\otimes}\,C^\infty_c(\mathbb{R}_+),\quad\mathcal{D}(G\times Y)'=\mathcal{S}(\mathbb{R})'\,\hat{\otimes}\,C^\infty_c(\mathbb{R}_+)'.
$$
The dual group is $\hat G=\mathbb{R}$. We choose $\Xi(\hat G)\doteq L^q(\mathbb{R})$ for some fixed $q\in[1,+\infty)$, and $\mathcal{Y}(Y)\doteq X(\mathbb{R}_+)$ a Banach space of distributions uniformly embedded in $C^\infty_c(\mathbb{R}_+)'$ in the sense of (\ref{UniformEmbed}).  The first important object is the space
$$
\Xi(\hat G,\mathcal{Y}(Y))\doteq\left\{\hat u\in L^1_\mathrm{loc}(\mathbb{R},X(\mathbb{R}_+))\;\vline\quad\|\varphi(\cdot)\|_{X(\mathbb{R}_+)}\in L^q(\mathbb{R})\right\}.
$$
We interpret $\hat u\in\Xi(\hat G,\mathcal{Y}(Y))$ as distributions $\hat u(\xi,y)$ on $\hat G\times Y=\Pi$, such that the map $\xi\mapsto\hat u(\xi,\cdot)$ is Bochner-measurable $\mathbb{R}\to X(\mathbb{R}_+)$, and the map $\xi\mapsto\|\hat u(\xi,\cdot)\|_{X(\mathbb{R}_+)}$ is from $L^q(\mathbb{R})$. It is easy to check that the (quasi-)lattice properties (\ref{LatticeXiNorm}) and (\ref{LatticeXiSpace}) are satisfied, therefore, by Lemma \ref{LatticeLemma}, we have that $\Xi(\hat G,\mathcal{Y}(Y))$ is a normed vector space with the norm
$$
\|\hat u\|_{\mathcal{Y},\Xi}=\|\|\hat u(\cdot)\|_{X(\mathbb{R}_+)}\|_q=\left(\int\limits_{-\infty}^\infty\|\hat u(\xi,\cdot)\|_{X(\mathbb{R}_+)}^q\right)^{\frac1q}.
$$
By Lemma \ref{EmbeddingLemma}, we have a continuous linear embedding
$$
\Xi(\hat G,\mathcal{Y}(Y))\hookrightarrow\mathcal{D}(\hat G\times Y)'=\mathcal{S}(\mathbb{R})'\,\hat{\otimes}\,C^\infty_c(\mathbb{R}_+)'.
$$
Since the condition (\ref{Xiu<|u|}) is trivially true, and $L^q(\mathbb{R})\hookrightarrow L^0(\mathbb{R})$ continuously, we have by Proposition \ref{XiMeasureClosed} that Cauchy sequences in $\Xi(\hat G,\mathcal{Y}(Y))$ converge at least in measure.

The mixed-Fourier-norm space $\mathcal{L}^{q;X}(\Pi)\doteq\mathcal{X}(G\times Y)$ is introduced as
$$
\mathcal{L}^{q;X}(\Pi)=\left\{u\in\mathcal{D}(G\times Y)'\;\vline\quad\hat u\in\Xi(\hat G,\mathcal{Y}(Y))\right\},
$$
with the norm $\|u\|=\|\hat u\|_{\mathcal{Y},\Xi}$. It follows immediately that
$$
\mathfrak{F}:\mathcal{L}^{q;X}(\Pi)\to\Xi(\hat G,\mathcal{Y}(Y))
$$
is an isometric isomorphism of normed vector spaces.

The function $\rho:\mathbb{R}\to(0,+\infty]$ is defined as
$$
\rho(\xi)=\left\|e^{-2\pi\langle\cdot,\xi\rangle}\right\|_{X(\mathbb{R}_+)},\quad\forall\xi\in\mathbb{R}.
$$
The weighted space $\Xi(\hat G,\rho)$ is given as
$$
\Xi(\hat G,\rho)=\left\{\hat u\in L^1_\mathrm{loc}(\mathbb{R})\;\vline\quad \hat u\rho\in L^q(\mathbb{R})\right\},
$$
with the norm $\|\hat u\|_{\Xi,\rho}=\|\hat u\rho\|_q$. By Lemma \ref{LatticeLemmaScalar}, $\Xi(\hat G,\rho)$ is a normed vector space (which was already clear at this point). Note, however, that the weight function $\rho^q$ need not be almost everywhere finite, in which case all elements of $\Xi(\hat G,\rho)$ must vanish on subsets where $\rho$ is essentially infinite. Since properties (\ref{Xi|u|<u}), (\ref{Xiu<|u|}), (\ref{u_|u|}) and (\ref{|u|_u}) are obviously true, if $\frac1\rho\in L^\infty_\mathrm{loc}(\mathbb{R})$, then by Proposition \ref{XirhoCompleteProp}, the normed vector space $\Xi(\hat G,\rho)$ is complete, i.e., a Banach space. For instance, for the weighted space $X(\mathbb{R}_+)=L^p(\mathbb{R}_+,\nu_\lambda)$ with $p\in[1,+\infty)$, the function $\rho$ is an explicitly computable function:
$$
\rho(\xi)=\begin{cases}
\xi^{-\frac{\lambda+1}{p}}\left(\frac{2^\lambda \Gamma (\lambda+2)}{p^{\lambda+1}}\right)^{\frac1{p}} & \mbox{for}\quad\xi>0,\\
+\infty & \mbox{for}\quad\xi\le0,
\end{cases}
$$
and $\frac1\rho\in C(\mathbb{R})$ (with the convention $\frac1{+\infty}=0$).

The mixed-Fourier-norm Bergman space $\mathcal{A}^{q;X}(\Pi)\doteq\mathcal{A}_\mathcal{X}(G\times Y)$ consists of all members of the space $\mathcal{L}^{q;X}(\Pi)$ which are holomorphic in $\Pi$,
$$
\mathcal{A}^{q;X}(\Pi)=\mathcal{L}^{q;X}(\Pi)\cap\mathrm{Hol}(\Pi)
$$
$$
=\left\{u\in\mathrm{Hol}(\Pi)\cap\mathcal{D}(G\times Y)'\;\vline\quad\hat u\in\Xi(\hat G,\mathcal{Y}(Y))\right\}.
$$
Suppose that the map
$$
\mathbb{R}\ni\xi\mapsto\frac{e^{-2\pi\langle\cdot,\xi\rangle}}{\rho(\xi)}\in X(\mathbb{R}_+)
$$
equals almost everywhere a Bochner-measurable function. This is clearly true, for instance, in the case $X(\mathbb{R}_+)=L^p(\mathbb{R}_+,\nu_\lambda)$ with $p\in[1,+\infty)$, where this map is actually continuous. By Proposition \ref{MainProp} we have that
$$
\mathfrak{F}:\mathcal{A}^{q;X}(\Pi)\to e^{-2\pi\langle\cdot,\cdot\rangle}\,\cdot\,\Xi(\hat G,\rho),\quad u(z)\mapsto e^{-2\pi\xi y}\,\hat u_0(\xi),
$$
is an isometric isomorphism of normed vector spaces. In particular, if $\Xi(\hat G,\rho)$ is complete (see above), then so is $\mathcal{A}^{q;X}(\Pi)$.

Let us finally turn to the Paley-Wiener-type properties. Suppose that the conditions of Proposition \ref{PWSuppProp} are met. This is obviously true for $X(\mathbb{R}_+)=L^p(\mathbb{R}_+,\nu_\lambda)$ with $p\in[1,+\infty)$, as we saw above. Then we have a Paley-Wiener-type support property,
\begin{equation}
\mathrm{supp}\,\hat u_0\subset[0,+\infty),\quad\forall u\in\mathcal{A}^{q;X}(\Pi).
\end{equation}
Moreover, suppose that the function
\begin{equation}
\frac{e^{-2\pi\langle y_0,\cdot\rangle}}{\rho}\in L^\infty(\mathbb{R})\label{Cond2}
\end{equation}
for some $y_0\ge0$. This holds for every $y_0>0$ in the case $X(\mathbb{R}_+)=L^p(\mathbb{R}_+,\nu_\lambda)$ with $p\in[1,+\infty)$. Then by Proposition \ref{PWBddProp} we have a Paley-Wiener-type boundedness property,
\begin{equation}
\sup_{y>y_1}\|u(\cdot,y)\|_\infty<\infty,\quad\forall u\in\mathcal{A}^{q;X}(\Pi),
\end{equation}
for every $y_1>y_0$.

\begin{remark}
This case of the mixed-Fourier-norm space associated with Cartesian coordinates on the half-plane is studied in more detail in the joint work of the authors of this article and Irina Smirnova, and will be published elsewhere. In particular, explicit formulas for equalities and representations of Paley-Wiener type and other results will be given. Here we present the basics of the definition of these spaces as an example of the application of our general approach and ideas.
\end{remark}

\subsection{Mixed-Fourier-norm space: the hyperbolic geometry case}\label{section-2}

The hyperbolic geometry will be realised on the domain $\Omega=\Pi\subset\mathbb{C}$ with the Lie group $G=\mathbb{R}$ acting freely as $w\mapsto e^xw$ on $w\in\Pi$. However, in this case our model space will be $\Gamma\doteq G\times Y=\mathbb{R}\times(0,\pi)\subset\mathbb{C}$, the (open) horizontal strip in the complex plane, with coordinates $z=x+\imath y\in\Gamma$, $x\in\mathbb{R}=G$, $y\in(0,2\pi)=Y$. We will use the bi-holomorphism $\Phi=\exp:\Gamma\to\Pi$,
$$
\Gamma\ni x+\imath y=z\longmapsto\exp(z)=e^{x+\imath y}=re^{\imath y}=w\in\Pi.
$$
We often use the polar coordinates $w=(r,y)\in\Pi$. Note that in this case both $y\to 0$ and $y\to2\pi$ correspond to $w\to\mathbb{R}$, while $x\to+\infty$ corresponds to $w\to+\imath\infty$. Thus, the $y$-direction is bounded, which renders Paley-Wiener-type properties irrelevant. This drawback will be reflected in somewhat limited scope of this model compared to the previous two.

On $Y=(0,\pi)$ we consider the measure
$$
d\nu_\lambda(y)=(\lambda+1)2^\lambda\sin^\lambda ydy
$$
for a fixed $\lambda>-1$. Our test function and distribution spaces will be
$$
\mathcal{D}(G\times Y)=\mathcal{S}(\mathbb{R})\,\hat{\otimes}\,C^\infty_c((0,\pi)),\quad\mathcal{D}(G\times Y)'=\mathcal{S}(\mathbb{R})'\,\hat{\otimes}\,C^\infty_c((0,\pi))'.
$$
The dual group is $\hat G=\mathbb{R}$. We choose $\Xi(\hat G)\doteq L^q(\mathbb{R})$ for some fixed $q\in[1,+\infty)$, and $\mathcal{Y}(Y)\doteq X((0,\pi))$ a Banach space of distributions uniformly embedded in $C^\infty_c((0,\pi))'$ in the sense of (\ref{UniformEmbed}).  The first important object is the space
$$
\Xi(\hat G,\mathcal{Y}(Y))\doteq\left\{\hat u\in L^1_\mathrm{loc}(\mathbb{R},X((0,\pi)))\;\vline\quad\|\varphi(\cdot)\|_{X((0,\pi))}\in L^q(\mathbb{R})\right\}.
$$
We interpret $\hat u\in\Xi(\hat G,\mathcal{Y}(Y))$ as distributions $\hat u(\xi,y)$ on $\hat G\times Y=\Gamma$, such that the map $\xi\mapsto\hat u(\xi,\cdot)$ is Bochner-measurable $\mathbb{R}\to X((0,\pi))$, and the map $\xi\mapsto\|\hat u(\xi,\cdot)\|_{X((0,\pi))}$ is from $L^q(\mathbb{R})$. It is easy to check that the (quasi-)lattice properties (\ref{LatticeXiNorm}) and (\ref{LatticeXiSpace}) are satisfied, therefore, by Lemma \ref{LatticeLemma}, we have that $\Xi(\hat G,\mathcal{Y}(Y))$ is a normed vector space with the norm
$$
\|\hat u\|_{\mathcal{Y},\Xi}=\|\|\hat u(\cdot)\|_{X((0,\pi))}\|_q=\left(\int\limits_{-\infty}^\infty\|\hat u(\xi,\cdot)\|_{X((0,\pi))}^q\right)^{\frac1q}.
$$
By Lemma \ref{EmbeddingLemma}, we have a continuous linear embedding
$$
\Xi(\hat G,\mathcal{Y}(Y))\hookrightarrow\mathcal{D}(\hat G\times Y)'=\mathcal{S}(\mathbb{R})'\,\hat{\otimes}\,C^\infty_c((0,\pi))'.
$$
Since the condition (\ref{Xiu<|u|}) is trivially true, and $L^q(\mathbb{R})\hookrightarrow L^0(\mathbb{R})$ continuously, we have by Proposition \ref{XiMeasureClosed} that Cauchy sequences in $\Xi(\hat G,\mathcal{Y}(Y))$ converge at least in measure.

The mixed-Fourier-norm space $\mathcal{X}(G\times Y)$ is introduced as
$$
\mathcal{X}(G\times Y)=\left\{u\in\mathcal{D}(G\times Y)'\;\vline\quad\hat u\in\Xi(\hat G,\mathcal{Y}(Y))\right\},
$$
with the norm $\|u\|=\|\hat u\|_{\mathcal{Y},\Xi}$. It follows immediately that
$$
\mathfrak{F}:\mathcal{X}(G\times Y)\to\Xi(\hat G,\mathcal{Y}(Y))
$$
is an isometric isomorphism of normed vector spaces.

The function $\rho:\mathbb{R}\to(0,+\infty]$ is defined as
$$
\rho(\xi)=\left\|e^{-2\pi\langle\cdot,\xi\rangle}\right\|_{X((0,\pi))},\quad\forall\xi\in\mathbb{R}.
$$
The weighted space $\Xi(\hat G,\rho)$ is given as
$$
\Xi(\hat G,\rho)=\left\{\hat u\in L^1_\mathrm{loc}(\mathbb{R})\;\vline\quad \hat u\rho\in L^q(\mathbb{R})\right\},
$$
with the norm $\|\hat u\|_{\Xi,\rho}=\|\hat u\rho\|_q$. By Lemma \ref{LatticeLemmaScalar}, $\Xi(\hat G,\rho)$ is a normed vector space (which was already clear at this point). Note, however, that the weight function $\rho^q$ need not be almost everywhere finite, in which case all elements of $\Xi(\hat G,\rho)$ must vanish on subsets where $\rho$ is essentially infinite. Since properties (\ref{Xi|u|<u}), (\ref{Xiu<|u|}), (\ref{u_|u|}) and (\ref{|u|_u}) are obviously true, if $\frac1\rho\in L^\infty_\mathrm{loc}(\mathbb{R})$, then by Proposition \ref{XirhoCompleteProp}, the normed vector space $\Xi(\hat G,\rho)$ is complete, i.e., a Banach space. For instance, for the weighted space $X((0,\pi))=L^p((0,\pi),\nu_\lambda)$ with $p\in[1,+\infty)$, the function $\rho$ is an explicitly computable function:
$$
\rho(\xi)=\left(\frac{\pi\Gamma(\lambda
+2)}{\left|\Gamma \left (\frac{\lambda
+2}{2}+\frac{p\xi i}{2}\right)\right|}\right)^\frac1{p}e^{-\frac{\pi\xi}{2}},\quad\forall\xi\in\mathbb{R},
$$
and $\frac1\rho\in C^\infty(\mathbb{R})$.

The holomorphic mixed-Fourier-norm space $\mathcal{A}_\mathcal{X}(G\times Y)$ consists of all members of the space $\mathcal{X}(G\times Y)$ which are holomorphic in $\Pi$,
$$
\mathcal{A}_\mathcal{X}(G\times Y)=\mathcal{X}(G\times Y)\cap\mathrm{Hol}(\Gamma)
$$
$$
=\left\{u\in\mathrm{Hol}(\Gamma)\cap\mathcal{D}(G\times Y)'\;\vline\quad\hat u\in\Xi(\hat G,\mathcal{Y}(Y))\right\}.
$$
Suppose that the map
$$
\mathbb{R}\ni\xi\mapsto\frac{e^{-2\pi\langle\cdot,\xi\rangle}}{\rho(\xi)}\in X((0,\pi))
$$
equals almost everywhere a Bochner-measurable function. This is clearly true, for instance, in the case $X(\mathbb{R}_+)=L^p((0,\pi),\nu_\lambda)$ with $p\in[1,+\infty)$, where this map is actually continuous. By Proposition \ref{MainProp} we have that
$$
\mathfrak{F}:\mathcal{A}_\mathcal{X}(G\times Y)\to e^{-2\pi\langle\cdot,\cdot\rangle}\,\cdot\,\Xi(\hat G,\rho),\quad u(z)\mapsto e^{-2\pi\xi y}\,\hat u_0(\xi),
$$
is an isometric isomorphism of normed vector spaces. In particular, if $\Xi(\hat G,\rho)$ is complete (see above), then so is $\mathcal{A}_\mathcal{X}(G\times Y)$.

Now let us translate all this back to the original domain $\Omega=\Pi$. We denote
$$
\Phi^*\mathcal{D}(\Pi)=\mathcal{D}(G\times Y),\quad\Phi^*\mathcal{D}(\Pi)'=\mathcal{D}(G\times Y)'.
$$
Distributions $u\in\mathcal{D}(\Pi)'$ are such that, formally speaking, $u(r,y)$ is tempered in the variable $\ln r$.

In order to make contact with the facts established above in the model space $\Gamma$, we will need to utilize the well-known relation between the Mellin transform over $\mathbb{R}_+$ and the Fourier transform over $\mathbb{R}$, and to translate the results to the language of Mellin transform $\mathfrak{M}:L^2(\mathbb{R}_+,\mu_\lambda)\to L^2(\mathbb{R})$,
$$
\mathfrak{M}f(\xi)=\int\limits_0^\infty r^{-2\pi\imath\xi+\frac\lambda2}f(r)dr,
$$
$$
d\mu_\lambda(r)=r^{\lambda+1}dr.
$$
As with Fourier transform, Mellin transform can be defined as $\mathfrak{M}:\mathcal{D}(\Pi)\to\mathcal{D}(G\times Y)$ by the extension of
$$
\mathfrak{M}\otimes\mathbf{1}:\exp^*\mathcal{S}(\mathbb{R})\,\hat{\otimes}\,C^\infty_c((0,\pi))\to\mathcal{S}(\mathbb{R})\,\hat{\otimes}\,C^\infty_c((0,\pi)),
$$
and as $\mathfrak{M}:\mathcal{D}(\Pi)'\to\mathcal{D}(G\times Y)'$ by duality. Explicitly,
$$
\tilde f(\xi,y)=\mathfrak{M}f(\xi,y)=\int\limits_0^\infty r^{-2\pi\imath\xi+\frac\lambda2}f(r,y)dr.
$$

With $\lambda>-1$ fixed before, let us introduce a linear bijection $\operatorname{T}_\lambda:C(\Pi)\to C(\Gamma)$ by setting
$$
\operatorname{T}_\lambda f(z)=e^{(\frac\lambda2+1)z}f(e^z),\quad\forall z\in\Gamma,\quad\forall f\in C(\Pi).
$$
Clearly, $\operatorname{T}_\lambda$ gives a bijection between holomorphic functions,
\begin{equation}
\operatorname{T}_\lambda:\mathrm{Hol}(\Pi)\to\mathrm{Hol}(\Gamma).\label{THolEq}
\end{equation}
Moreover, one can check directly that
$$
\left\|\operatorname{T}_\lambda f(\cdot,y)\right\|_p=\left\|f(\cdot,y)\right\|_{L^p(\mathbb{R}_+,\mu_{\lambda,p})},
$$
$$
d\mu_{\lambda,p}(r)=r^{(\frac\lambda2+1)p-1}dr,
$$
$$
\forall y\in(0,\pi),\quad\forall f\in C_c(\Pi),\quad\forall p\in[1,\infty).
$$
The importance of $\operatorname{T}_\lambda$ lies in the fact that it connects Mellin transform with Fourier transform. More precisely, we have the following property, which can be checked directly:
\begin{equation}
\mathfrak{M}f(\xi,y)=e^{-\imath y(1+\frac\lambda2)}\mathfrak{F}\operatorname{T}_\lambda f(\xi,y),\label{MellinFourier}
\end{equation}
$$
\forall(\xi,y)\in\Gamma,\quad\forall f\in C_c(\Gamma).
$$
It can be checked directly that the following commutative diagram is true in the sense of continuous embeddings of topological vector spaces,
$$
\begin{tikzcd}
C^\infty_c(\Pi) \arrow[d,"{\operatorname{T}_\lambda}"] \arrow[r,hook] & \mathcal{D}(\Pi) \arrow[d,"{\operatorname{T}_\lambda}"] \arrow[r,hook] & L^2(\Pi,\mu_\lambda) \arrow[d,"{\operatorname{T}_\lambda}"] \arrow[r,hook] & \mathcal{D}(\Pi)' \arrow[d,"{\operatorname{T}_\lambda}"] \arrow[r,hook] & C^\infty_c(\Pi)' \arrow[d,"{\operatorname{T}_\lambda}"]\\
C^\infty_c(\Gamma) \arrow[r,hook] & \mathcal{D}(G\times Y) \arrow[r,hook] & L^2(\Gamma) \arrow[r,hook] & \mathcal{D}(G\times Y)' \arrow[r,hook] & C^\infty_c(\Gamma)'
\end{tikzcd}
$$
In the above diagram, $\operatorname{T}_\lambda$ is everywhere an isomorphism. From the property (\ref{MellinFourier}) it follows that
$$
\mathfrak{M}\,\mathcal{D}(\Pi)=\mathcal{D}(G\times Y),
$$
\begin{equation}
\mathfrak{M}\,\mathcal{D}(\Pi)'=\mathcal{D}(G\times Y)'.\label{MellinDistEq}
\end{equation}
If $u\in\mathcal{D}(G\times Y)'$ and $f\in\mathcal{D}(\Pi)'$ such that $u=\operatorname{T}_\lambda f$, then $\hat u=\tilde f\in\mathcal{D}(G\times Y)$.

The mixed-Fourier-norm space $\mathcal{L}^{q;X}(\Pi)$ is defined to be the preimage of $\mathcal{X}(G\times Y)$,
$$
\operatorname{T}_\lambda\mathcal{L}^{q;X}(\Pi)=\mathcal{X}(G\times Y).
$$
More precisely,
$$
\mathcal{L}^{q;X}(\Pi)=\left\{f\in\mathcal{D}(\Pi)'\;\vline\quad\tilde f\in L^1_\mathrm{loc}(\mathbb{R},X((0,\pi)),\quad\|f\|<\infty\right\},
$$
with the norm
$$
\|f\|=\left(\int\limits_{-\infty}^\infty\|\tilde f(\xi)\|_{X((0,\pi))}^q\right)^{\frac1q}.
$$

Finally, the mixed-Fourier-norm Bergman space $\mathcal{A}^{q;X}(\Pi)$ is defined as the preimage of the holomorphic mixed-Fourier-norm space $\mathcal{A}_\mathcal{X}(G\times Y)$,
$$
\operatorname{T}_\lambda\mathcal{A}^{q;X}(\Pi)=\mathcal{A}_\mathcal{X}(G\times Y).
$$
In precise terms,
$$
\mathcal{A}^{q;X}(\Pi)=\left\{f\in\mathrm{Hol}(\Pi)\cap\mathcal{D}(\Pi)'\;\vline\quad\tilde f\in\Xi(\hat G,X((0,\pi)))\right\}.
$$

If the space $\mathcal{A}_\mathcal{X}(G\times Y)$ is complete (see above), that implies the comleteness of the mixed-Fourier-norm Bargman space $\mathcal{A}^{q;X}(\Pi)$. We restate again, that Mellin transform $f\mapsto\tilde f$ is an isometric isomorphism of Banach spaces
$$
\mathcal{A}^{q;X}(\Pi)\to e^{-2\pi\langle\cdot,\cdot\rangle}\,\cdot\,\Xi(\hat G,\rho),\quad f(z)\mapsto e^{-2\pi\xi y}\,\tilde f_0(\xi).
$$

\begin{remark}
The present case the mixed-Fourier-norm space associated with polar coordinates on the half-plane is studied in more detail in the work of Irina Smirnova, and will be published separately. In particular, explicit formulas for equalities and representations of Paley-Wiener type and other results will be given. Here we present the basics of the definition of these spaces as an example of the application of our general approach and ideas.
\end{remark}

\section*{Acknowledgements}
The work is supported by the Ministry of Education and Science of Russia, agreement No. 075-02-2024-1427. The first author is supported by the FWO Odysseus 1 grant G.0H94.18N: Analysis and Partial Differential Equations and by the Methusalem programme of the Ghent University Special Research Fund (BOF) (Grant number 01M01021).

\section*{Data availability statement}
The authors confirm that all data generated or analysed during this study are included in this article.

\end{document}